\documentclass{siamart1116}

\usepackage{graphicx,epstopdf}
\usepackage{subcaption}
\usepackage{multirow}
\usepackage[colorinlistoftodos]{todonotes}
\usepackage{marginnote}
\usepackage{algorithmic} 
\Crefname{ALC@unique}{Line}{Lines} 
\newcommand{\REMOVE}[1]{}

\usetikzlibrary{calc}

\providecommand{\mytodo}[2][]{{%
\let\marginpar\marginnote
\reversemarginpar
\renewcommand{\baselinestretch}{0.8}%
\todo[#1]{#2}}}
\definecolor{greeo}{RGB}{31,113,9}
\newcommand{\rstumin}[1]{{\color{greeo}{RST: #1}}}

\usepackage{lipsum}
\usepackage{amsfonts}
\usepackage{graphicx}
\usepackage{epstopdf}
\usepackage{algorithmic}
\ifpdf
  \DeclareGraphicsExtensions{.eps,.pdf,.png,.jpg}
\else
  \DeclareGraphicsExtensions{.eps}
\fi

\numberwithin{theorem}{section}

\newcommand{\TheTitle}{A Monolithic Algebraic Multigrid Framework for Multiphysics Applications with examples from Resistive {MHD}}
\newcommand{\TheAuthors}{Ohm, Wiesner, Cyr, Hu, Shadid, and Tuminaro}

\headers{Monolithic AMG Framework for Multiphysics Systems}{\TheAuthors}

\title{{A Monolithic Algebraic Multigrid Framework for Multiphysics Applications with examples from 
Resistive {MHD}}\thanks{Submitted to the editors 3/2021.
\funding{This work was supported by the U.S.~Department of Energy, Office of Science, Office of Advanced Scientific
       Computing Research, Applied Mathematics program and by the U.S. Department of Energy, Office of Science, 
       Office of Advanced Scientific Computing Research and Office of Fusion Energy Sciences, Scientific Discovery through Advanced Computing (SciDAC) program.
       Sandia National Laboratories is a multimission laboratory managed and operated by National Technology
       and Engineering Solutions of Sandia, LLC., a wholly owned subsidiary of Honeywell International, Inc.,
       for the U.S. Department of Energy's National Nuclear Security Administration under grant~DE-NA-0003525.
       This paper describes objective technical results and analysis.  Any subjective views or opinions that
       might be expressed in the paper do not necessarily represent the views of the U.S. Department of Energy
       or the United States Government.}}}

\author{
  Peter Ohm\thanks{Sandia National Laboratories, P.O. Box 5800, MS 1320, Albuquerque, NM 87185
    (\email{pohm@sandia.gov}, \email{eccyr@sandia.gov}, \email{jnshadi@sandia.gov}).} 
  \and
  Tobias A. Wiesner\thanks{Leica Geosystems AG, Heinrich-Wild-Strasse 201, 9435 Heerbrugg/SG, Switzerland (\email{tobias@tawiesn.de})}
  \and
  Eric C. Cyr\footnotemark[2]
  \and
  Jonathan J. Hu\thanks{Sandia National Laboratories, P.O. Box 969, MS 9159, Livermore, CA 94661
    (\email{jhu@sandia.gov}, \email{rstumin@sandia.gov}).}
  \and
  John N. Shadid{$^\dagger$}\thanks{Department of Mathematics and Statistics, University of New Mexico, Albuquerque NM, 87123}
  \and
  Raymond S. Tuminaro\footnotemark[4]
}

\usepackage{amsopn}

\ifpdf
\hypersetup{
  pdftitle={\TheTitle},
  pdfauthor={\TheAuthors}
}
\fi

\begin{document}

\maketitle

\begin{abstract}
  A multigrid framework is described for multiphysics applications. 
  The framework allows one to construct, adapt, and tailor a monolithic
  multigrid methodology  to different linear systems coming from discretized 
  partial differential equations. The main idea centers on developing multigrid
  components in a blocked fashion where each block corresponds to separate
  sets of physical unknowns and equations within the larger
  discretization matrix. Once defined, these components are ultimately
  assembled into a monolithic multigrid solver for the entire system. 
  We demonstrate the potential of the framework by applying it to 
  representative linear solution sub-problems arising from resistive MHD.
\end{abstract}



\providecommand{\vels}{u}      
\providecommand{\vel}{{\bf \vels}} 
\providecommand{\pre}{p}       
\providecommand{\visc}{\mu}    
\providecommand{\dens}{\rho} 
\providecommand{\resistivity}{\eta} 
\providecommand{\magnperm}{\mu_0} 
\providecommand{\viscstress}{{\bf \Pi}} 
\providecommand{\totalstress}{{\bf \mathcal{T}}} 
\providecommand{\fluidstress}{{\bf T}} 
\providecommand{\magnstress}{{\bf T}_M} 
\providecommand{\magns}{B} 
\providecommand{\magn}{{\bf \magns}} 
\providecommand{\current}{{\bf J}} 
\providecommand{\inductionflux}{{\bf \mathcal{F}}} 
\providecommand{\lagr}{\psi} 
\providecommand{\res}{\text{r}} 
\providecommand{\resv}{{\bf r}} 
\providecommand{\lundquist}{S} 

\providecommand{\bgsdamping}{\omega}

\providecommand{\solv}{\mathbf{x}}	
\providecommand{\solvO}{\solv_{\textnormal{0}}}
\providecommand{\solvl}{\solv_{\textnormal{1}}}

\providecommand{\resvO}{\resv_{\textnormal{0}}}
\providecommand{\resvl}{\resv_{\textnormal{1}}}

\providecommand{\rhsv}{\mathbf{b}}
\providecommand{\rhsvO}{\mathbf{b}_{\textnormal{0}}}
\providecommand{\rhsvl}{\mathbf{b}_{\textnormal{1}}}

\providecommand{\level}{\ell}

\section{Introduction \& Motivation}

Multigrid methods are among the fastest most scalable techniques for 
solving the sparse linear systems that arise from the discretization
of many different systems of partial differential equations (PDEs)~\cite{TrOoSc2001}. 
The basic multigrid idea is to accelerate convergence to a linear system's solution by performing 
relaxation (i.e., simple iterative techniques) on a hierarchy of
different resolution systems.  Algebaic multigrid methods (AMG)
are popular as they build the hierarchy automatically requiring
little effort from the application developer.  While algebraic multigrid 
methods have been successfully applied to many complex problems, further
developments are needed to more robustly adapt them to multiphysics
PDE systems. In fact, the original AMG research was focused on 
scalar PDEs, and most of the theory relies on
the linear system being symmetric positive definite.

Multigrid's rapid convergence relies on constructing its 
components such that relaxation sweeps applied to the different fidelity
versions are complementary to each other.  That is, errors not easily damped 
by relaxation on one version can be damped effectively by relaxation
on another version. 
However, constructing AMG components with desirable complementary
properties can be quite complicated for PDE systems when the 
coupling between different physical unknowns is significant.
AMG challenges arise in designing effective relaxation procedures
as well as in designing the grid transfer operators used to project
the high fidelity version of the matrix problem to lower resolutions.
A discussion of AMG issues for PDE systems and some different
approaches can be found in ~\cite{FuSt02}. 
While there are important AMG themes when considering
PDE systems, the {\it best} AMG choices are typically
application dependent. For this reason, it is critical 
that an AMG framework can support the implementation and exploration
of a wide range of AMG adaptations for complex multiphysics
PDE systems. In this paper, we describe such a framework that has
been implemented in the {\sf MueLu} package~\cite{MueLu,MueLuURL}, found within the 
Trilinos libraries~\cite{trilinos-article}. We demonstrate
its utility in the context of solving difficult linear systems associated 
with magnetohydrodynamics (MHD) simulations.  In particular, we highlight how the 
framework allows one to adapt the method to different scenarios. In one
case, a specialized block relaxation method is devised that has significant
computational advantages over a more black-box relaxation technique.
In another case, we show how to adapt the solver to address situations
where different finite element basis functions are used to represent
the different physical fields of the MHD system. 

The resistive magnetohydrodynamics (MHD) model provides a base-level continuum description of 
the dynamics of conducting fluids in the presence of electromagnetic fields. This model 
is useful in the context of plasma physics systems coming from both 
naturally occurring dynamics (e.g. astrophysics and planetary dynamos) 
and man-made systems (e.g. magnetic confinement fusion) \cite{GoedbloedPoedts2004}.  
The governing partial differential equations for the resistive  MHD model consist
of conservation of mass, momentum and energy augmented by the low-frequency Maxwell's equations.
These systems are strongly coupled, highly nonlinear and characterized by coupled physical phenomena that span a
very large range of length- and time-scales.  These characteristics make the scalable, robust, 
accurate, and efficient computational solution of these systems extremely challenging. 
From this point of view, fully-implicit formulations, coupled with effective robust nonlinear/linear iterative solution methods,
become attractive, as they have the potential to provide stable,
higher-order time-integration of these complex multiphysics systems. These methods can
follow the dynamical time-scales of interest, as opposed to
time-scales determined by fast normal modes that 
severely impact numerical stability, but do not control temporal order-of-accuracy.
The existence of fast normal modes in MHD make the associated algebraic systems very stiff \cite{GoedbloedPoedts2004}, and thus
difficult to solve iteratively using modern iterative techniques \cite{chacon-pop-08-3dmhd,shadid-jcp-10_mhd,JardinReview2011,shadid2016scalable}. Despite these challenges, much effort has been invested in
the last decade towards the development of more fully-implicit, nonlinearly coupled solution methods for MHD, with the goal of
enhancing robustness, accuracy, and efficiency 
(see e.g. \cite{keppens-ijnm-99-imhd,toth-astroph-98-imhd,chacon-JCP-rmhd,luisMHD,
ShumlakLoverich_2003,chacon-cpc-04-mhd_discret,knoll-prl-06-ic,chacon-pop-08-3dmhd,shadid-jcp-10_mhd,Shumlaketal_multifluid_2011,JardinReview2011,shadid2016scalable,CyrABF2DResistiveMHD2010}).
To efficiently address the unique linear systems challenge associated with MHD simulations, special
solver adaptations must be considered to leverage the underlying character of the equations.
The framework described in this paper is intended to facilitate a range of different solver adaptations
that might be appropriate for different MHD scenarios.

In Section \ref{sec:mhd} we introduce the MHD equations and the accompanying discrete systems of equations.
In Section \ref{section:amg} we give an overview of the algebraic multigrid method.
We discuss the implementation of these algebraic multigrid methods to multiphysics PDE systems
in a truly monolithic manner through the use of blocked operators in Section \ref{sec:blockAMG}.
We demonstrate the numerical and computational performance benefits of this approach
on various test problems and present the results in Section \ref{sec:experiments}.
Finally we end with concluding remarks in Section \ref{sec:conclusion}.

\section{MHD equations}
\label{sec:mhd}


The model of interest for this paper is the 3D resistive iso-thermal MHD equations including dissipative terms for the 
momentum and magnetic induction equations  \cite{GoedbloedPoedts2004}.  
The system of equations is:
%

\begin{align}
\frac{\partial \bigl(\dens \vel\bigr)}{\partial t} + \nabla \cdot \bigl[ \dens \vel \otimes \vel - \totalstress \bigr] &= {\bf 0},
\label{Eq: MomentumEq} \\
\frac{\partial {\dens }} {\partial t} + \nabla \cdot \bigl[ \dens {\vel } \bigr]  &= { 0},
\label{Eq: ContinuityEq} \\
%
%
\frac{\partial\magn}{\partial t} + \nabla \cdot \Bigl[ \vel \otimes \magn - { \bf B} \otimes \vel - \frac{\resistivity}{\magnperm} \bigl(\nabla \magn - (\nabla \magn)^T \bigr)  \Bigr] &= {\bf 0}.
\label{Eq: MagEvolutionB}
\end{align}
Here $\vel$ is the plasma center of mass velocity; $\dens$ is the mass
density; $\totalstress$ is the (total) stress tensor;    $\magn$ is the magnetic induction (here after also termed the magnetic field) that is subject to the 
divergence-free involution $\nabla \cdot \magn = 0$. 
The associated current, $\current$, is obtained from Amp\`{e}re's law as $\current = \frac{1}{\mu_0} \nabla \times \magn$.
In its simplest form, the resistive MHD equations are completed with definition of the total stress tensor, $\totalstress = \fluidstress + \magnstress$
composed of the  fluid and magnetic stress tensors, $\fluidstress$, $\magnstress$ that are given by,
%

%
%
%
%
%
%
%

\begin{subequations}
\begin{align}
\textbf{T} &=  -[ \pre{\bf I} -\frac{2}{3} \mu (\nabla \cdot \vel)]  {\bf I} + \visc[\nabla \vel +  \nabla \vel^{T}],\\
\magnstress &=  \frac{1}{\magnperm} \magn \otimes  \magn -    \frac{1}{2 \magnperm} \| \magn \|^2  {\bf I}.
\end{align}
\end{subequations}
Here, $\pre$ is the plasma pressure 
and ${\bf I}$ is the identity tensor. 
The physical parameters in this model are  the plasma viscosity, $\visc$,  the resistivity, $\resistivity$, and  the
magnetic permeability of free space, $\magnperm$. Finally, for convenience of notation,   a tensor induction flux, $\inductionflux$, that is defined as 
\begin{equation}
 \inductionflux = \vel \otimes \magn - \magn \otimes \vel - \frac{\resistivity}{\magnperm} \left(\nabla \magn - (\nabla \magn)^T  \right)  + \psi {\bf I}
\label{Eq: TensorInductionFlux}
\end{equation}
is introduced.

Satisfying the solenoidal involution $\nabla \cdot\magn = 0$ in the discrete representation to machine precision 
is a topic of considerable  interest in both structured and
unstructured finite-volume and unstructured finite-element contexts
(see e.g. \cite{toth-jcp-00-divB,Dender02,chacon-cpc-04-mhd_discret}).
In the formulation discussed in this study, a scalar Lagrange multiplier ($\psi$ in Eq.~\eqref{Eq: TensorInductionFlux})
is introduced into the induction equation that enforces the 
solenoidal involution as a divergence free constraint on the magnetic field. 
This procedure is common in both finite volume (see e.g. \cite{Dender02,toth-jcp-00-divB,chacon-cpc-04-mhd_discret})
and finite element methods ({see e.g. \cite{Codina_MHD2006,Codina_MHD2011,BadiaCodinaPlanas13,shadid2016scalable}). 

This paper focuses on the incompressible limit of this system, i.e., $\nabla \cdot \vel = 0$. 
This limit is useful in 
the modeling of various applications such as
low-Lundquist-number liquid-metal MHD flows \cite{Moreau90,Davidson01}, and high-Lundquist-number, large-guide-field 
fusion plasmas \cite{strauss-76-phfl-rmhd,hazeltine-85-phfl-rmhd,drake-phfl-84-rmhd}.
This limit is characteristic of low flow-Mach-number applications for compressible systems as well, and is the most challenging 
algorithmically because of the presence of the elliptic incompressibility constraint. 
However, the 
stabilized FE formulation that is presented, the strongly-coupled Newton-Krylov nonlinear iterative solvers, and the 
fully-coupled algebraic multilevel preconditioners also work in the variable density low-Mach-number compressible case.
Together the incompressibility constraint with the solenoid 
involution, enforced as a constraint, produces a dual saddle point structure for the systems of equations \cite{shadid2016scalable}.

The  3D resistive MHD equations in residual form with the introduction of the scalar Lagrange multiplier and the incompressibility assumption are given by

\begin{subequations}
\label{tab:governing-transport-pdes}
\begin{align}
\resv_m &=   \frac{\partial  \dens  \vel}{\partial t} +    \nabla  \cdot \bigl[ \dens  \vel \otimes \vel - \totalstress\bigr]  = {\bf 0} & \textnormal{Momentum}\\
\res_P &=   \dens  \nabla \cdot \vel =  0 & \textnormal{Continuity Constraint}\\
\resv_I &= \frac{\partial \magn}{\partial t} + \nabla \cdot \inductionflux = {\bf 0} & \textnormal{Magnetic Induction}\\
\res_{\lagr} &= \nabla \cdot \magn  = 0 & \textnormal{Solenoidal Constraint}
\end{align}
\end{subequations}
The spatial discretization of the governing equations is based on the VMS FE method \cite{TJRH95,DoneaHuerta2002}.
The semi-discretized system is integrated in time with a method of lines approach based on BDF schemes. 
The weak form of the VMS / Stabilized FE formulation for the resistive MHD equation (see Eq.~\eqref{tab:governing-transport-pdes})
is 
given by
\begin{subequations}
\label{tab:governing-weak-transport-pdes}
{\small
\begin{align}
\ {\bf F}^h_{{\bf u}} &= 
\!\! {\displaystyle \int_{\Omega}} \hspace*{.1cm}  \!\! {\bf w}^h \cdot {\resv}^h_{m} d\Omega
+\!\! {\displaystyle \sum_e} {\displaystyle \int_{\Omega_e}} \hspace*{.1cm} 
\!\! \rho \hat \tau_{m}  {\resv}^h_{m}  \otimes {\bf u}^h : \nabla {\bf w}^hd \Omega
+\!\! {\displaystyle \sum_e} {\displaystyle \int_{\Omega_e}} \hspace*{.1cm} 
\!\! \hat \tau_P ( \nabla \cdot {\bf w}^h) {{\res}^h_P} ~d \Omega, \\
F^h_P & = 
\!\! {\displaystyle \int_{\Omega}} \hspace*{.1cm} \!\! q^h {\res}^h_P d\Omega
+\!\! {\displaystyle \sum_e} {\displaystyle \int_{\Omega_e}} \hspace*{.1cm}
\!\! \rho \hat \tau_m \nabla q^h \cdot {\resv}^h_m d \Omega,\\
{\bf F}^h_{{I}} &= 
\!\! {\displaystyle \int_{\Omega}} \hspace*{.1cm}  \!\! {{\bf C}^h \cdot {\resv}^h_{I}} d\Omega
-\!\! {\displaystyle \sum_e} {\displaystyle \int_{\Omega_e}} \hspace*{.1cm} 
\!\! \hat \tau_{I} ( {{\resv}^h_{I}  \otimes {\bf u}^h} - {\bf u}^h \otimes  {{\resv}^h_{I}} ) : \nabla {\bf C}^h~ d \Omega
+\!\! {\displaystyle \sum_e} {\displaystyle \int_{\Omega_e}} \hspace*{.1cm} 
\!\! \hat \tau_{\bf \psi} ( \nabla \cdot {\bf C}^h){{\res}^h_{\psi}} d \Omega, \label{eq:vms-maginduc} \\
F^h_{\psi} &= 
\!\! {\displaystyle \int_{\Omega}} \hspace*{.1cm} \!\! {s^h  {\res}^h_{\psi} }d\Omega
+\!\! {\displaystyle \sum_e} {\displaystyle \int_{\Omega_e}} \hspace*{.1cm}
\!\! \hat \tau_{ I} \nabla s^h \cdot {\resv}^h_{I} d \Omega ,
\end{align}
}
\end{subequations}
\noindent where $\hat \tau_i$ are the stabilization parameters.
Here $[{\bf w}^h,q^h,{\bf C}^h,s^h]$ are the FE  weighting functions for the velocity, pressure, magnetic field and the Lagrange multiplier respectively. The sum ${\small \sum_e}$ indicates the integrals are taken only over element interiors ${\small \Omega_e}$ and integration by parts is not performed. A full development and examination of this formulation is presented in~\cite{shadid2016scalable}.

\subsection{Brief Overview of Discrete Systems of Equations}

To provide context to the solution methods that follow later, 
we present a brief discussion of the equation structure.
Here, the focus is on the structure of VMS terms generated by 
the induction equation and in the enforcement of the solenoidal constraint through the Lagrange multiplier, $\lagr$. 
Specifically, 
the VMS 
weak form of the  solenoidal constraint
in expanded form, is detailed in \cite{shadid2016scalable} as
\begin{equation}
F_{\lagr} = 
\displaystyle{ \int_{\Omega}
s  \bigl(\nabla \cdot \magn\bigr) ~d\Omega~+ } 
  \displaystyle{ \sum_e \int_{\Omega_e}
\hat \tau_I \nabla s \cdot \biggl( \frac{\partial \magn}{\partial t} + 
\nabla\cdot \mathcal{F} \biggr) d \Omega}.
\label{eqn:ExpandedContinuityWeakForm}
\end{equation} 
%
This expression includes a weak Laplacian operator acting on the Lagrange multiplier
\begin{equation}
{ {L_{\lagr}} =  {\displaystyle \sum_e} {\displaystyle \int_{\Omega_e}}
\hat \tau_I  \nabla s \cdot  \nabla \lagr  d \Omega }.
\label{eqn:PressureWeakLaplacian}
\end{equation}
This term is the analogue of the weak pressure Laplacian 
${ \footnotesize {L_{P} =  {\sum_e} {\int_{\Omega_e}} \hspace*{.1cm} \dens \hat \tau_m  \nabla \Phi \cdot  \nabla \pre  d \Omega }}$ appearing in 
the total mass continuity equation (see the general discussion for stabilized FE CFD in \cite{DoneaHuerta2002} and our previous development in \cite{shadid-jcp-10_mhd,shadid2016scalable} in the context of 
2D MHD). These VMS operators are critical in the elimination of oscillatory modes from the null space of the resistive MHD saddle point system for both $(\vel,\pre)$ and 
$(\magn,\lagr)$ and allow equal-order interpolation of all the unknowns (see \cite{DoneaHuerta2002} for incompressible CFD and \cite{Codina_MHD2006,Codina_MHD2011,BadiaCodinaPlanas13} for resistive MHD and  a coercivity proof of stability). 

The final term of the stabilized form of the magnetics induction equation~\eqref{eq:vms-maginduc} involves a weak divergence-type operator (substitute $\res_\psi$).
This term adds to the stability of the VMS form \cite{Codina_MHD2006,Codina_MHD2011,BadiaCodinaPlanas13,shadid2016scalable} and also enhances the ability to iteratively invert the magnetic induction sub-block in the Jacobain matrix by physics-based and approximate block factorization methods \cite{EHSST06,CyrTeko2016}.

A finite element (FE) discretization of the stabilized equations
gives rise to a system of coupled, nonlinear, non-symmetric algebraic
equations, the numerical solution of which can be very
challenging. These equations are linearized using an inexact form of
Newton's method leading to a block system with the following form
%
\providecommand{\discgrad}{{\bf G^T}}
\providecommand{\lap}{{\bf L}}
\providecommand{\nsu}{{\bf J_\vel}}
\providecommand{\nsgrad}{{\bf G}} 
\providecommand{\nsd}{{\bf D}}
\providecommand{\nslap}{\bf \lap_P}

\providecommand{\magb}{{\bf J_I}}
\providecommand{\maglap}{{\bf \lap_\lagr}}

\providecommand{\coupY}{{\bf Y}}
\providecommand{\coupZ}{{\bf Z}}


\begin{equation}
\begin{bmatrix}
\nsu & \nsgrad & \coupZ & {\bf 0} \\
\nsd & \nslap  & {\bf 0} & {\bf 0} \\
\coupY & {\bf 0} & \magb & \nsgrad \\
{\bf 0} & {\bf 0} & \discgrad & \maglap
\end{bmatrix} 
\left[
\begin{array} {c}
{\bf  \delta \hat \vel} \\
{\delta \hat \pre} \\
{\bf \delta \hat B} \\
{\bf \delta \hat \psi}
\end{array}
\right]
=
-\left[
\begin{array} {c}
\resv_\vel \\
\res_\pre \\
\resv_I \\
\res_\lagr \\
\end{array}
\right].
\label{eqn:4x4BlockMatrix}
\end{equation}

The block matrix, $\nsu$, corresponds to the discrete transient, convection, diffusion and stress terms acting on the unknowns ${\bf \delta \hat \vel}$; the matrix,
$\nsgrad$, corresponds to the discrete gradient operator; $\nsd$, the discrete representation of the continuity equation terms with velocity (note for
a true incompressible flow this would be the divergence operator denoted $\discgrad$); the block matrix, $\magb$, corresponds to the discrete transient, convection, diffusion terms acting on magnetic induction, and the matrices, $\nslap, \maglap$, are the stabilization Laplacian's mentioned above.  The right hand side vectors  contain the residuals for Newton's method. We solve for the solution increments $\delta\hat\vel$ and $\delta\hat\pre$ for the velocity and pressure as well as for $\bf \delta\hat B$ and $\delta\hat\psi$ which represent the magnetics field and the Lagrange multipliers.
The $\nslap$ and $\maglap$ operators
help facilitate the solution of the linear systems with a number of algebraic and domain
decomposition type preconditioners that rely on non-pivoting ILU factorization, Jacobi relaxation or Gauss-Seidel as
sub-domain solvers \cite{shadid1999,shadid2016scalable,shadid2004}. 

The difficulty of producing robust and efficient preconditioners to \eqref{eqn:4x4BlockMatrix}
has motivated many different types of
decoupled solution methods.  Often, transient schemes
combine semi-implicit methods with fractional-step
(operator splitting) approaches or use fully-decoupled solution
strategies \cite{toth-astroph-98-imhd,keppens-ijnm-99-imhd,aydemir-jcp-85-imhd,m3d,jardin-pop-05-mhd,harned-kerner-si,schnack-barnes-si,harned-mikic-si,nimrod2,hujeirat-mnras-98-irmhd,AlegraALE08}. In these cases, the motivation is to reduce
memory usage and to produce a simplified equation set
for which efficient solution strategies already exist. Unfortunately, 
these simplifications place significant limitations
on the broad applicability of these methods.
A detailed presentation of the characteristics of
different linear and nonlinear solution strategies is beyond our current scope.
Here, we wish to highlight that our approach 
of fully-coupling the resistive MHD PDEs in the nonlinear solver preserves the inherently
strong coupling of the physics with the goal to produce a more robust
solution methodology \cite{shadid1999,shadid-jcp-10_mhd,shadid2016scalable}. Preservation of this strong coupling, however,
places a significant burden on the linear solution procedure. 


\section{Multigrid methods}

\providecommand{\MphysicsA}{gray!70!white}
\providecommand{\MphysicsB}{gray!50!white}
\providecommand{\MphysicsC}{gray!30!white}
\providecommand{\MphysicsCouplingA}{lightgray!50!white}
\providecommand{\MphysicsCouplingC}{lightgray!40!white}
\providecommand{\MphysicsCouplingB}{lightgray!30!white}

\providecommand{\A}{\mathbf{A}}
\providecommand{\AOO}{\A_{00}}
\providecommand{\AOl}{\A_{01}}
\providecommand{\AlO}{\A_{10}}
\providecommand{\All}{\A_{11}}

\providecommand{\Prol}{\textnormal{\textbf{P}}}
\providecommand{\POO}{\Prol_{00}}
\providecommand{\Pll}{\Prol_{11}}

\providecommand{\Restr}{\textnormal{\textbf{R}}}
\providecommand{\ROO}{\Restr_{00}}
\providecommand{\Rll}{\Restr_{11}}

%
%
Multigrid methods are based on the fact that many simple iterative methods are effective at eliminating high-frequency error components relative to the ``mesh'' resolution used for discretization.  
The basic multigrid idea is to introduce a hierarchy of discrete approximations to the PDE problem 
employing different resolution meshes. The newly introduced coarser approximations are used
to accelerate the solution process on the finest mesh as lower-frequency error components on the finest
mesh will appear to be high-frequency (relative to the grid resolution) on a particular mesh within 
the hierarchy.  That is, different error components are efficiently reduced by essentially applying a simple iterative method to the appropriate resolution approximation.

Multigrid methods generally come in two varieties: geometric multigrid (GMG) and algebraic multigrid (AMG). Typically with GMG , applications supply a hierarchy of meshes, discrete PDE operators, interpolation operators to
transfer solutions from coarser resolutions to finer ones,
and restriction operators to transfer finer resolution residuals to coarser
levels.  In geometric multigrid, inter-grid transfers are based on
geometric relationships, such as using linear interpolation to 
define a finer level approximation from a coarse solution.
With AMG methods, coarse level information and inter-grid transfers are
developed automatically by analyzing the supplied fine level discretization
matrix. This normally involves a combination of graph heuristics to coarsen
followed by some approximation algorithm to develop inter-grid transfers 
with the aim of more accurately transferring information associated with small
magnitude eigenvalues of the discrete operator,  which is often low frequency.
This paper focuses on AMG as it can be more easily adapted to complex
application domains by non-multigrid scientists.

\begin{figure}
\caption{Multigrid V-cycle pseudo-code and graphical representation for a 3-level method \label{MG algorithm}}
\begin{tabbing}

\= $\textsf{MGV}(
u_{\ell}, b_{\ell},\ell):$ \\[2pt]
\>\hskip .15in\= {\it if} $ \ell \ne \ell_{max}  $ \\
\>   \>\hskip .15in\=$u_{\ell} $ \hskip .11in\=
        $\leftarrow \mathcal{S}_\ell^\textnormal{pre} ({\bf A}_{\ell},u_{\ell}, b_{\ell}) $ \\
\>\>\> $r_\ell$ \> $\leftarrow b_\ell - {\bf A}_\ell u_\ell $ \\
\>\>\> $u_{\ell+1}$ \> $\leftarrow 0 $ \\
\>\>\> $u_{\ell+1}$ \> $\leftarrow \textsf{MGV}(
      u_{\ell+1}, \Restr_{\ell\rightarrow{\ell+1}} r_{\ell}, 
                      \ell\hskip-3pt+\hskip-3pt1)$ \\
\>\>\> $u_{\ell}$ \> $\leftarrow u_{\ell} + \Prol_{\ell+1\rightarrow{\ell}} 
u_{\ell+1} $\\
\>\>\> $u_{\ell}$ \> $\leftarrow \mathcal{S}_\ell^\textnormal{post} ({\bf A}_{\ell},u_{\ell}, b_{\ell})  $ \\[-2pt]
\>\> {\it else} \\[-3pt]
\>\>\> $ u_{\ell} \leftarrow {\bf A}_{\ell}^{-1} b_{\ell} $\\[-20pt]
\end{tabbing}
\vskip -1.6in \hskip 2.05in \begin{tikzpicture}
    %
	%

	\begin{scope}[shift={(8,1.3)},scale=0.2] 
     	\tikzstyle{restrict}=[line width=0.4mm,->]	
     	\tikzstyle{prolongate}=[line width=0.4mm,->]	     	
      \begin{scope}[shift={(0,0)},scale=0.4]
        \begin{scope}[every node/.append style={yslant=-0.5,xslant=1},yslant=-0.5,xslant=1]
        	\draw[line width=0.5mm,color=black,fill=\MphysicsA] (0,0) rectangle (7,-7);
        	\node[color=white] at (3.5,-3.5) {\scalebox{0.6}{$\A^{(2)}$}};        	        	
		\node (left-paren) at (-0.6,-3.5) {$\left[\vphantom{\rule{0.6cm}{0.46cm}}\right.$};
		\node (right-paren) at (7.6,-3.5) {$\left.\vphantom{\rule{0.6cm}{0.46cm}}\right]$};
		\end{scope}
      \end{scope}
      
      \draw[restrict,->] (-5,3) -- (0,-1.5); 
      \draw[prolongate,->] (0.1,-1.5) -- (5,3); 

      \begin{scope}[shift={(-5,5)},scale=0.6,opacity=0.9]
        \begin{scope}[every node/.append style={yslant=-0.5,xslant=1},yslant=-0.5,xslant=1]        
        	\draw[line width=0.5mm,color=black,fill=\MphysicsA] (0,0) rectangle (7,-7);
        	\node[color=white] at (3.5,-3.5) {$\A^{(1)}$};        	
		\node (left-paren) at (-0.6,-3.3) {$\left[\vphantom{\rule{0.7cm}{0.6cm}}\right.$};
		\node (right-paren) at (7.6,-3.3) {$\left.\vphantom{\rule{0.7cm}{0.6cm}}\right]$};
		\end{scope}
      \end{scope}
      \begin{scope}[shift={(5,5)},scale=0.6,opacity=0.9]
        \begin{scope}[every node/.append style={yslant=-0.5,xslant=1},yslant=-0.5,xslant=1]
        	\draw[line width=0.5mm,color=black,fill=\MphysicsA] (0,0) rectangle (7,-7);
        	\node[color=white] at (3.5,-3.5) {$\A^{(1)}$};        	
		\node (left-paren) at (-0.6,-3.3) {$\left[\vphantom{\rule{1cm}{0.6cm}}\right.$};
		\node (right-paren) at (7.6,-3.3) {$\left.\vphantom{\rule{1cm}{0.6cm}}\right]$};		
		\end{scope}
      \end{scope}

      \draw[restrict,->] (-10,7.2) -- (-5,3); 
      \draw[prolongate,->] (5,3) -- (10,7.2); 
                  
     \begin{scope}[shift={(-10,10)},scale=0.8,opacity=0.9]  
        \begin{scope}[every node/.append style={yslant=-0.5,xslant=1},yslant=-0.5,xslant=1]
        	\draw[line width=0.5mm,color=black,fill=\MphysicsA] (0,0) rectangle (7,-7);
        	\node[color=white] at (3.5,-3.5) {\scalebox{1.4}{$\A^{(0)}$}};
		\node (left-paren) at (-0.5,-3.5) {$\left[\vphantom{\rule{1cm}{0.9cm}}\right.$};
		\node (right-paren) at (7.5,-3.5) {$\left.\vphantom{\rule{1cm}{0.9cm}}\right]$};
		\end{scope}
      \end{scope}
      \begin{scope}[shift={(10,10)},scale=0.8,opacity=0.9]
        \begin{scope}[every node/.append style={yslant=-0.5,xslant=1},yslant=-0.5,xslant=1]
        	\draw[line width=0.5mm,color=black,fill=\MphysicsA] (0,0) rectangle (7,-7);
        	\node[color=white] at (3.5,-3.5) {\scalebox{1.4}{$\A^{(0)}$}};
		\node (left-paren) at (-0.5,-3.5) {$\left[\vphantom{\rule{1cm}{0.9cm}}\right.$};
		\node (right-paren) at (7.5,-3.5) {$\left.\vphantom{\rule{1cm}{0.9cm}}\right]$};
		\end{scope}
	  \end{scope}
	  
      \begin{scope}[every node/.append style={yslant=-0.5,xslant=1},yslant=-0.5,xslant=1]
      \node at (-12.0,-2.5) {$\mathcal{S}_{0}^\textnormal{pre}$};
      \node at (-5.5,-4.0) {$\mathcal{S}_{1}^\textnormal{pre}$};
      \node at (2.4,-5) {$\mathcal{S}_{2}^\textnormal{pre/post}$};
      \node at (+5.0,4.5) {$\mathcal{S}_{1}^\textnormal{post}$};      
      \node at (+4.0,12) {$\mathcal{S}_{0}^\textnormal{post}$};      
      \end{scope}  
      
      \begin{scope}[every node/.append style={yslant=0.5,xslant=0},yslant=0,xslant=0]
            \node at (-12.3,2.2) {$\Restr_{0\rightarrow 1}$};
            \node at (-6.6,-1.2) {$\Restr_{1\rightarrow 2}$};
	  \end{scope}  
	  
	  \begin{scope}[every node/.append style={yslant=-0.5,xslant=0.0},yslant=0,xslant=0]
            \node at (6.8,-1.6) {$\Prol_{2\rightarrow 1}$};
            \node at (13.1,1.8) {$\Prol_{1\rightarrow 0}$};
	  \end{scope}  
      
	\end{scope}
	\end{tikzpicture}%
\end{figure}

Figure~\ref{MG algorithm} depicts what is referred to as a multigrid
V-cycle to solve a linear system ${\bf A}_{\ell} u_{\ell} = b_{\ell}$.
Subscripts distinguish between different resolutions.
$\Prol_{\ell+1\rightarrow\ell}$
interpolates from level $\ell\hskip -.03in+\hskip -.03in1$ to
level $\ell$.
$\Restr_{\ell\rightarrow\ell+1}$ restricts from level $\ell$ to
level $\ell\hskip -.03in+\hskip -.03in1$.
${\bf A}_{\ell}$ is the discrete problem on level $\ell$ and for 
coarse levels is defined by a Petrov-Galerkin projection
$$
{\bf A}_{\ell+1}
= 
\Restr_{\ell\rightarrow\ell+1} {\bf A}_{\ell} \Prol_{\ell+1\rightarrow\ell} .
$$
{$\mathcal{S}_\ell^\textnormal{pre}$}
and
{$\mathcal{S}_\ell^\textnormal{post}$}
denote a basic iterative scheme
(e.g., Gauss-Seidel) that is applied to damp or relax some error components.
The overall efficiency is governed by the interplay of the two main multigrid ingredients: inter-grid transfer operators and the smoothing methods.
One can either use several sweeps with the multigrid V-cycle as a standalone solver, or one can apply a V-cycle sweep as a preconditioner within an outer Krylov based linear iterative solver.  
For a general overview on multigrid methods the reader is 
referred to~~\cite{BrHeMc2000,TrOoSc2001}
and the references therein.

For PDE and multiphysics systems, applying AMG to the entire PDE system
(i.e., monolithic multigrid)
can be problematic, especially when the coupling between different 
solution types (e.g., pressures and velocities) is strong.
Classical simple iterative methods may not necessarily reduce all 
oscillatory error components and might even amplify some oscillatory
error components. In fact, methods requiring the inversion
of the matrix diagonal (e.g., the Jacobi iteration) are not even well 
defined when applied to incompressible fluid formulations that give rise 
to zeros on the matrix diagonal. Further, many standard AMG algorithms for
defining inter-grid transfers might lead to transfers that
do not accurately preserve smooth functions. For example, methods such
as smoothed aggregation rely on a Jacobi-like step to generate smooth 
inter-grid transfer basis functions. This Jacobi step, however, will
obviously not generate smoother basis functions when the Jacobi iteration 
is not well defined (or when the matrix diagonal is {\it small} in a 
relative sense). 
In even worse situations, the AMG software may not be applicable to the PDE 
system when different FE basis functions are used to represent different 
fields within the system. This is because most AMG codes employ a simple 
technique to address PDE systems. This technique relies on the different 
equations to be effectively discretized in a similar way. 
In particular, AMG coarsening algorithms (or for our approach aggregation 
algorithms) are often applied by first grouping all DoFs at each mesh node 
together and then applying the coarsening algorithms to the graph induced
from the block matrix. An advantage to this approach is that all unknowns
at mesh points are coarsened in identically the same fashion. However,
the approach cannot be applied when the number of DoFs associated with
different fields varies or if all DoFs are not co-located, e.g., 
when using quadratic basis functions to represent velocities while 
pressures employ only linear basis functions. In previous versions 
of our AMG software, the only alternative to this simple PDE system
technique would be to completely ignore the multiphysics coupling and 
effectively treat the entire system as if it is a scalar PDE, which 
almost always leads to disastrous convergence rates.


\label{section:amg}

\section{Truly monolithic block multigrid for the MHD equations}
\label{sec:blockAMG}

In this section we propose a multigrid method for multiphysics
systems that can be represented by block matrices allowing one to
adapt both the relaxation algorithms and the grid transfer construction
algorithms to the structure of the system.

\subsection{Block matrices and multigrid for PDE systems}
As illustrated in the previous MHD discussion, 
PDE systems are often represented by block matrices as in, for example, equation~\eqref{eqn:4x4BlockMatrix}.
More generally, block systems can be written as
\begin{equation}
\begin{bmatrix}
A_{00} & A_{01} & \cdots & A_{0N} \\
A_{10} & A_{11} & \cdots & A_{1N} \\
\vdots & \vdots & \ddots & \vdots \\ 
A_{N0} & A_{N1} & \cdots & A_{NN} \\
\end{bmatrix}
\begin{bmatrix}
x_0 \\
x_1 \\
\vdots \\
x_N \\
\end{bmatrix}
=
\begin{bmatrix}
b_0 \\
b_1 \\
\vdots \\
b_N \\
\end{bmatrix}
\label{eqn:block-prec-form}
\end{equation}
where each component of the vector $x_i$ is a field in the multiphysics PDE
and the sub-matrices $A_{ij}$ are approximations to operators in the governing 
equations.

One preconditioning approach to PDE systems follows a so-called {\it physics-based}
strategy (see Figure~\ref{fig:mgmphmodelsA2}). These techniques can be viewed as approximate block factorizations 
(involving Schur complement approximations) to the block matrix equations \eqref{eqn:block-prec-form}.
The factorizations are usually constructed based on the underlying physics.
Here, different AMG V cycle sweeps are used to approximate the different sub-matrix inverses that appear within the 
approximate block factors.
As the sub-matrices correspond to single physics or scalar PDE operators, 
application-specific modifications to the multigrid algorithm are often not necessary. This makes the physics-based strategy 
particularly easy to implement as one can leverage ready-to-use multigrid packages.
Several methods that follow this popular strategy are described in ~\cite{ElSiWa2005} and the 
references therein.
While the physics-based approach has some practical advantages, the efficiency of 
the preconditioner relies heavily on how well the coupling between different equations within the PDE
is approximated by the block factorization.

In this paper, we instead consider a monolithic multigrid alternative to a physics-based 
approach.  A monolithic scheme applies a multigrid algorithm to the entire block PDE system and so the AMG scheme effectively develops a hierarchy of 
block PDE matrices associated with different resolutions. Figure~\ref{fig:mgmphmodels2} graphically illustrates the two contrasting
approaches highlighting the key potential advantage to a monolithic approach. 
\begin{figure}[htbp]
\centering
\caption{Two multigrid approaches to address multiphysics applications.}
\label{fig:mgmphmodels2}
\begin{subfigure}[t]{0.45\textwidth}\centering
\begin{tikzpicture}[scale=0.4,every node/.style={minimum size=1cm}] 
\begin{scope}[xshift=0]

\draw[black,thick,dotted] (1,2.5) -- (-0.5,-3);    
\draw[black,thick,dotted] (-5,2.5) -- (-3.5,-3);   
\draw[black,thick,dotted] (-2,4) -- (-2,-2.2);     

  \draw[black,thick,dotted] (3,3.5) -- (3,1.5);    
  \draw[black,thick,dotted] (1.1,2.5) -- (1.5,0.75);    
  \draw[black,thick,dotted] (5,2.5) -- (4.5,0.75);    

    \begin{scope}[xshift=-30,yshift=-120,scale=0.5,
    	    every node/.append style={
    	    yslant=0.5,xslant=-1},yslant=0.5,xslant=-1
    	             ]
		\fill[\MphysicsA,fill opacity=.9] (0,2) rectangle (3,5);
        \draw[black,thick] (0,2) rectangle (3,5);
        \draw (1.5,3.5) node[scale=0.5]{\tiny$\AOO^{(3)}$};
    \end{scope}  
    
    \begin{scope}[xshift=-20,yshift=-80,scale=0.666,
    	    every node/.append style={
    	    yslant=0.5,xslant=-1},yslant=0.5,xslant=-1]
		\fill[\MphysicsA,fill opacity=.9] (0,2) rectangle (3,5);
        \draw[black,thick] (0,2) rectangle (3,5);
        \draw (1.5,3.5) node[scale=0.75]{\tiny$\AOO^{(2)}$};
    \end{scope}  
    
    \begin{scope}[xshift=-10,yshift=-40,scale=0.833,
    	    every node/.append style={
    	    yslant=0.5,xslant=-1},yslant=0.5,xslant=-1]
		\fill[\MphysicsA,fill opacity=.9] (0,2) rectangle (3,5);
        \draw[black,thick] (0,2) rectangle (3,5);
        \draw (1.5,3.5) node[scale=1.0]{\tiny$\AOO^{(1)}$};
    \end{scope}

    \begin{scope}[xshift=26,yshift=-30,scale=0.7,
    	    every node/.append style={
    	    yslant=0.5,xslant=-1},yslant=0.5,xslant=-1]
		\fill[\MphysicsB,fill opacity=.9] (3,0) rectangle (5,2);
        \draw[black,thick] (3,0) rectangle (5,2);
        \draw (4,1) node[scale=0.5]{\tiny$\All^{(1)}$};
    \end{scope}

    \begin{scope}[xshift=0,yshift=0,scale=1,
    	    every node/.append style={
    	    yslant=0.5,xslant=-1},yslant=0.5,xslant=-1]
		\fill[\MphysicsA,fill opacity=.9] (0,2) rectangle (3,5);
		\fill[\MphysicsB,fill opacity=.9] (3,0) rectangle (5,2);	
		\fill[\MphysicsCouplingA,fill opacity=.9] (0,0) rectangle (3,2);
		\fill[\MphysicsCouplingA,fill opacity=.9] (3,2) rectangle (5,5);
        \draw[black,thick] (0,2) rectangle (3,5);
        \draw[black,thick] (3,0) rectangle (5,2);	
        \draw[black,thick] (0,0) rectangle (3,2);
        \draw[black,thick] (3,2) rectangle (5,5);
        \draw (1.5,3.5) node{$\AOO$};
        \draw (1.5,1) node{$\AOl$};
        \draw (4,3.5) node{$\AlO$};
        \draw (4,1) node{$\All$}; 
    \end{scope}     
  \draw[black,thick,dotted] (-2,1) -- (-2,-3.7);   
  \draw[black,thick,dotted] (3,1.5) -- (3,0);   
  \end{scope}    
\end{tikzpicture}
\caption{Physics-based approach with multiple multigrid approximations to sub-matrix inverses within an approximate block factorization.  Example with 4 and 2 multigrid levels for $\AOO$ and $\All$.}
\label{fig:mgmphmodelsA2}
\end{subfigure}
~~~\begin{subfigure}[t]{0.45\textwidth}\centering
\begin{tikzpicture}[scale=0.4,every node/.style={minimum size=1cm},on grid]
\begin{scope}[]

\draw[black,thick,dotted] (5,2.5) -- (2.5,-3);
\draw[black,thick,dotted] (-5,2.5) -- (-2.5,-3);
\draw[black,thick,dotted] (0,5) -- (0,-1.6);

    \begin{scope}[xshift=0,yshift=-120,scale=0.5,
    	    every node/.append style={
    	    yslant=0.5,xslant=-1},yslant=0.5,xslant=-1
    	             ]
		\fill[\MphysicsA,fill opacity=.9] (0,2) rectangle (3,5);
		\fill[\MphysicsB,fill opacity=.9] (3,0) rectangle (5,2);	
		\fill[\MphysicsCouplingA,fill opacity=.9] (0,0) rectangle (3,2);
		\fill[\MphysicsCouplingA,fill opacity=.9] (3,2) rectangle (5,5);
        \draw[black,thick] (0,2) rectangle (3,5);
        \draw[black,thick] (3,0) rectangle (5,2);	
        \draw[black,thick] (0,0) rectangle (3,2);
        \draw[black,thick] (3,2) rectangle (5,5);
        \draw (1.5,3.5) node[scale=0.5]{\tiny$\AOO^{(3)}$};
        \draw (1.5,1) node[scale=0.5]{\tiny$\AlO^{(3)}$};
        \draw (4,3.5) node[scale=0.5]{\tiny$\AOl^{(3)}$};
        \draw (4,1) node[scale=0.5]{\tiny$\All^{(3)}$}; 
    \end{scope}  
    
    \begin{scope}[xshift=0,yshift=-80,scale=0.666,
    	    every node/.append style={
    	    yslant=0.5,xslant=-1},yslant=0.5,xslant=-1
    	             ]
		\fill[\MphysicsA,fill opacity=.9] (0,2) rectangle (3,5);
		\fill[\MphysicsB,fill opacity=.9] (3,0) rectangle (5,2);	
		\fill[\MphysicsCouplingA,fill opacity=.9] (0,0) rectangle (3,2);
		\fill[\MphysicsCouplingA,fill opacity=.9] (3,2) rectangle (5,5);
        \draw[black,thick] (0,2) rectangle (3,5);
        \draw[black,thick] (3,0) rectangle (5,2);	
        \draw[black,thick] (0,0) rectangle (3,2);
        \draw[black,thick] (3,2) rectangle (5,5);
        \draw (1.5,3.5) node[scale=0.75]{\tiny$\AOO^{(2)}$};
        \draw (1.5,1) node[scale=0.75]{\tiny$\AlO^{(2)}$};
        \draw (4,3.5) node[scale=0.75]{\tiny$\AOl^{(2)}$};
        \draw (4,1) node[scale=0.75]{\tiny$\All^{(2)}$}; 
    \end{scope}  
    
    \begin{scope}[xshift=0,yshift=-40,scale=0.833,
    	    every node/.append style={
    	    yslant=0.5,xslant=-1},yslant=0.5,xslant=-1
    	             ]
		\fill[\MphysicsA,fill opacity=.9] (0,2) rectangle (3,5);
		\fill[\MphysicsB,fill opacity=.9] (3,0) rectangle (5,2);	
		\fill[\MphysicsCouplingA,fill opacity=.9] (0,0) rectangle (3,2);
		\fill[\MphysicsCouplingA,fill opacity=.9] (3,2) rectangle (5,5);
        \draw[black,thick] (0,2) rectangle (3,5);
        \draw[black,thick] (3,0) rectangle (5,2);	
        \draw[black,thick] (0,0) rectangle (3,2);
        \draw[black,thick] (3,2) rectangle (5,5);
        \draw (1.5,3.5) node{\tiny$\AOO^{(1)}$};
        \draw (1.5,1) node{\tiny$\AlO^{(1)}$};
        \draw (4,3.5) node{\tiny$\AOl^{(1)}$};
        \draw (4,1) node{\tiny$\All^{(1)}$}; 
    \end{scope}

    \begin{scope}[xshift=0,yshift=0,scale=1,
    	    every node/.append style={
    	    yslant=0.5,xslant=-1},yslant=0.5,xslant=-1
    	             ]
		\fill[\MphysicsA,fill opacity=.9] (0,2) rectangle (3,5);
		\fill[\MphysicsB,fill opacity=.9] (3,0) rectangle (5,2);	
		\fill[\MphysicsCouplingA,fill opacity=.9] (0,0) rectangle (3,2);
		\fill[\MphysicsCouplingA,fill opacity=.9] (3,2) rectangle (5,5);
        \draw[black,thick] (0,2) rectangle (3,5);
        \draw[black,thick] (3,0) rectangle (5,2);	
        \draw[black,thick] (0,0) rectangle (3,2);
        \draw[black,thick] (3,2) rectangle (5,5);
        \draw (1.5,3.5) node{$\AOO$};
        \draw (1.5,1) node{$\AlO$};
        \draw (4,3.5) node{$\AOl$};
        \draw (4,1) node{$\All$}; 
    \end{scope}     
  \draw[black,very thick,dashed] (0,0) -- (0,-4.2);    
  \end{scope}    
\end{tikzpicture}
\caption{Monolithic multigrid with PDE coupling represented on all multigrid levels. }
\label{fig:mgmphmodelsB2}
\end{subfigure}
\end{figure}
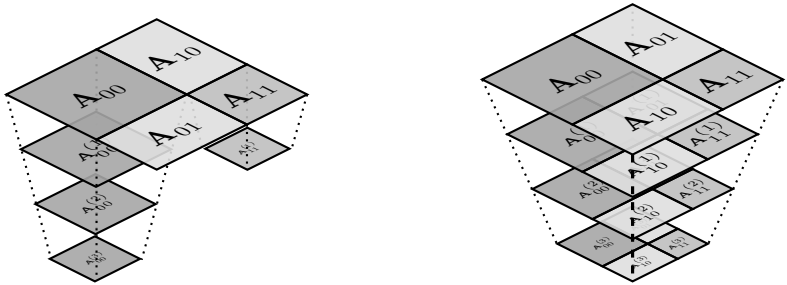
In particular, the 
cross-coupling defined by $\AOl$ and $\AlO$ is explicitly represented on all multigrid hierarchy levels 
with a monolithic approach as opposed to the physics-based scheme that attempts to represent coupling 
within an approximate Schur complement (whose inverse might then be approximated via multigrid sweeps). 
That is, a physics-based approach relies heavily on being able to develop effective Schur complements,
which can be non-trivial for complex applications. On the other hand, a monolithic approach 
introduces its own set of application-specific mathematical challenges such as the construction of 
relaxation procedures for monolithic systems and the development of coarsening schemes for different 
fields within a multiphysics system.  The design of efficient multiphysics preconditioners 
often requires one to make use of the specific knowledge about the block structure, the mathematical models of the underlying physics, and the problem-specific coupling of the equations.  The mathematical
challenges of multiphysics systems are often further compounded by non-trivial software challenges.
Unfortunately, most AMG packages cannot be customized
to particular multi-physics scenarios without having an in-depth knowledge of the AMG software. It is for this reason that many prefer a physics-based approach.

In the next sub-sections, we propose a monolithic algorithm for the MHD equations that can be
easily adopted and customized using the MueLu package within the Trilinos framework~\cite{MueLu,MueLuURL}.
Though we focus on a concrete MHD case, we cannot emphasize enough the importance of the software's generality in facilitating a monolithic approach for those with limited
knowledge of the multigrid package internals.

\subsection{Algebraic representation of the MHD problem}


In adapting a monolithic multigrid strategy to the MHD equations, a natural approach would be to interpret the system \eqref{eqn:4x4BlockMatrix} as a $2\times 2$ block system where the Navier-Stokes equations are separated from the Maxwell equations. That is, we treat the Navier-Stokes part and the Maxwell part as separate entities that are coupled by the off-diagonal blocks as shown in the $2\times 2$ block representation of Figure \ref{fig:2x2blocksystem}. In this way, we can leverage existing ideas/solvers for 
the Navier-Stokes equations and for the Maxwell equations.
\begin{figure}
\caption{Representation of MHD system \eqref{eqn:4x4BlockMatrix} as a $2\times 2$ block system}
\label{fig:2x2blocksystem}
\begin{center}
   	\begin{tikzpicture}
        \begin{scope}[shift={(6,0)},scale=0.7, every node/.append style={yslant=0,xslant=0},yslant=0,xslant=0]
   	   	\tikzstyle{style_fillAA}=[line width=0.2mm,color=black,fill=\MphysicsA]
   	   	\tikzstyle{style_fillAB}=[line width=0.2mm,color=black,fill=\MphysicsCouplingA]
	   	\tikzstyle{style_fillBA}=[line width=0.2mm,color=black,fill=\MphysicsCouplingA]   	   	
	   	\tikzstyle{style_fillBB}=[line width=0.2mm,color=black,fill=\MphysicsB]
        \draw[style_fillAA] (0,0) rectangle (2,-2);
		\draw[style_fillBB] (2.2,-2.2) rectangle (4.2,-4.2);

		\draw[style_fillAB] (2.2,0) rectangle (4.2,-2);			
		\draw[style_fillBA] (0,-2.2) rectangle (2,-4.2);			

		\draw[color=black,dashed] (0,-1.2) -- (2,-1.2);
		\draw[color=black,dashed] (2.2,-1.2) -- (4.2,-1.2);
		\draw[color=black,dashed] (1.2,0) -- (1.2,-2);
		\draw[color=black,dashed] (1.2,-2.2) -- (1.2,-4.2);		
		\draw[color=black,dashed] (0,-3.4) -- (2,-3.4);
		\draw[color=black,dashed] (2.2,-3.4) -- (4.2,-3.4);
		\draw[color=black,dashed] (3.4,0) -- (3.4,-2);
		\draw[color=black,dashed] (3.4,-2.2) -- (3.4,-4.2);

        \node at (0.6,-0.6)[text depth=0.45ex,text height=1.35ex] {$\nsu$}; \node at (1.6,-0.6)[text depth=0.45ex,text height=1.35ex] {$\nsgrad$};
        \node at (0.6,-1.6)[text depth=0.45ex,text height=1.35ex] {$\nsd$}; \node at (1.6,-1.6)[text depth=0.45ex,text height=1.35ex] {$\nslap$};

        \node at (2.8,-2.8)[text depth=0.45ex,text height=1.35ex] {$\magb$};     \node at (3.8,-2.8)[text depth=0.45ex,text height=1.35ex] {$\nsgrad$};
        \node at (2.8,-3.8)[text depth=0.45ex,text height=1.35ex] {$\discgrad$}; \node at (3.8,-3.8)[text depth=0.45ex,text height=1.35ex] {$\maglap$};        
        
        \node at (2.8,-0.6)[text depth=0.45ex,text height=1.35ex] {$\coupZ$};     
        \node at (0.6,-2.8)[text depth=0.45ex,text height=1.35ex] {$\coupY$};
        
		\node (left-paren) at (-0.25,-2.1) {$\left[\vphantom{\rule{0.2cm}{1.6cm}}\right.$};
		\node (right-paren) at (4.45,-2.1) {$\left.\vphantom{\rule{0.2cm}{1.6cm}}\right]$};
		\node (rightvec)[anchor=west] at (4.7,-2.1) {$\begin{bmatrix}
		{\bf  \delta \hat \vel} \\
    {\delta \hat \pre} \\
    {\bf \delta \hat B} \\
    {\delta \hat \psi}\end{bmatrix}=-
    \begin{bmatrix}
    \resv_\vel \\
    \res_\pre \\
    \resv_I \\
    \res_\lagr \\
    \end{bmatrix}$};	
		\end{scope}
	\end{tikzpicture}%
\end{center}
\end{figure}
In making this $2 \times 2$ decomposition, we are effectively emphasizing the 
significance of the coupling between fields within the Navier-Stokes block and 
within the Maxwell block as compared to the coupling between Navier-Stokes
unknowns and Maxwell unknowns. This is due to the importance of the coupling
constraint equations (e.g., incompressibility conditions involving velocities
or contact constraints \cite{wiesner2015}) to the associated evolution equations. 
These constraint equations often give rise to saddle-point like block systems.
Of course, there are physical situations where the coupling between the 
Navier-Stokes equations and the Maxwell equations is quite significant and so
a block $2 \times 2 $ decomposition might be less appropriate. There might also
be situations where coupling relationships are more complex. 
If, for example, the MHD equations are embedded in another larger more complex 
multiphysics problem, then one might need to consider a hierarchy of of coupling
configurations, which might require different arrangements/blocking of multigrid 
ingredients. Given the problem specific nature of multiphysics preconditioning,
our emphasis here is on the importance of a flexible software framework
to facilitating different types of blocking within the preconditioner.

For the remainder of the paper the block notation
\begin{equation} \label{eq:2x2system}
\begin{bmatrix} \AOO & \AOl \\ \AlO  & \All \end{bmatrix}\begin{bmatrix}\solvO\\\solvl\end{bmatrix}=
-
    \begin{bmatrix}
    \mathbf{b}_0 \\ \mathbf{b}_1
    \end{bmatrix},
\end{equation}
is used representing the corresponding blocks from Figure \ref{fig:2x2blocksystem}. The velocity and pressure increments ${\bf  \delta \hat \vel}$ and ${\bf  \delta \hat \pre}$ are grouped in $ \solvO$ and the Maxwell information is represented by $\solvl$, respectively. In a similar way, the block notation for the residual vector is adopted.


\subsection{Monolithic multigrid ingredients for volume-coupled problems}

The multiphysics solver that we propose is generally applicable to 
{\it volume-coupled} problems.  Volume coupled means that the different physics 
blocks are defined on the same domain. Volume-coupled examples include, e.g, thermo-structure-interaction (TSI) problems (see \cite{danowski2013}) or in our case the MHD 
equations. 
Specifically, within our MHD formulation all physics equations (i.e., Navier-Stokes and Maxwell parts) are defined throughout the entire domain.
This is in contrast to
interface-coupled problems such as fluid-structure interaction (FSI) applications (see \cite{gee2011,Deparis2016c,Langer2016a,Jodlbauer2019a}) or structural contact problems (see \cite{wiesner2015,Wiesner_ContactSP}) where different equation sets are valid over distinct domains that are only coupled through a common interface. From a multigrid perspective special interface coarsening methods are necessary for interface-coupled applications, which is not the focus of this paper.

Two multigrid ingredients must be specified to fully define the monolithic solver:
the inter-grid transfer operators and the relaxation or smoother  procedures. 

\subsubsection{Inter-grid transfers for the MHD system}
Following the $2\times 2$ decomposition of the MHD system from equation~\eqref{eq:2x2system}, we consider 
rectangular block diagonal inter-grid transfer operators
\begin{equation}
	\Restr_{i\rightarrow i+1}=\begin{bmatrix} \ROO &  \\   & \Rll \end{bmatrix} \textnormal{  and  }	\Prol_{i+1\rightarrow i}=\begin{bmatrix} \POO &  \\   & \Pll \end{bmatrix}
	\label{eq:blocktransfers}
\end{equation}
for restriction and prolongation between multigrid levels respectively.
The basic idea uses the {\sf MueLu} multigrid package to produce
grid transfers for the Navier-Stokes equations and the Maxwell
equations and then leverages {\sf MueLu}'s flexibility to combine these
grid transfer operators into a composite block diagonal operator.

The block perspective allows us to use completely separate invocations
of the multigrid package to build individual components (e.g., $\POO$ and $\Pll$)
and then combine or compose them together.  
As discussed earlier, we rely on underlying core AMG kernels that 
are applicable to either a scalar PDE or a PDE system with co-located unknowns.
Mixed finite element  schemes may not normally satisfy this co-located
requirement,  so the ability to separately invoke these core
components to produce $\POO$ and $\Pll$ alleviates this restriction, 
allowing us to apply monolithic AMG to a wider class of PDE systems.
That is, a mixed basis function discretization can be approached without
having to erroneously treat the entire system as a scalar PDE.
We demonstrate this capability at the end of
Section \ref{sec:experiments} through a mixed formulation utilizing
Q2/Q2 VMS for the Navier-Stokes degrees of freedom and
Q1/Q1 VMS for the Maxwell degrees of freedom.

In the case where the unknowns between blocks are co-located 
(e.g., Q1/Q1 VMS for both Navier-Stokes and Maxwell) we have the option to
correlate the grid transfer construction by having the multigrid invocations share the same
aggregates (or coarsening definition) as depicted in Figure \ref{fig:mgmphaggs}.
In this way, we guarantee that there is a one-to-one relationship of the coarse
Maxwell degrees of freedoms and the associated Navier-Stokes degrees of freedom.
That is, we obtain the same coarsening rate for both, and so the ratio
between Navier-Stokes and Maxwell degrees of freedom is constant on all multigrid  levels.
It should be noted that this is relatively straight-forward when all equations
are defined on the same mesh using a first-order nodal finite element 
discretization method.  Thus, there are 8 degrees of freedom at each mesh node 
(four associated with fluid flow and four associated with electromagnetics).

Within our aggregation approach, we first build standard aggregates using the graph of 
sub-block $\AOO$ on level $i$ as input.
Next we build separate aggregates using the graph of sub-block $\All$ on level $i$ as input.
However, if the degrees of freedom are co-located, instead of building new aggregates we have the option to
reconstruct the corresponding aggregates for the Maxwell equations by cloning the Navier-Stokes aggregation
information as depicted in Figure~\ref{fig:mgmphaggs}.
Specifically, mesh vertices on a given level are
assigned to \textit{aggregates} $\mathcal{A}_\ell^i$ such that
\begin{displaymath}
  \bigcup_{i=1}^{N_{\ell+1}} \mathcal{A}_\ell^i = \left\{ 1,...,N_\ell  \right\} \ , \ 
  \mathcal{A}_\ell^i \cap \mathcal{A}_\ell^j = \emptyset \ , \ 1 \leq i < j \leq N_{\ell+1} \ ,
\end{displaymath}
where $N_\ell$ denotes the number of mesh vertices on level $\ell$.
Each aggregate $\mathcal{A}_\ell^i$ on level $\ell$ gives rise to one node on level $\ell\hskip -.02in +\hskip -.02in 1$.  The $\mathcal{A}_\ell^i$ are formed by applying
greedy algorithms to the graph associated with a matrix discretization.
Typically, one wants aggregates to be approximately the same size
and roughly spherical in shape (for isotropic problems).
Piecewise-constant interpolation can then be defined over each aggregate for each 
solution component.

\newcommand{\convexpath}[2]{
[   
    create hullnodes/.code={
        \global\edef\namelist{#1}
        \foreach [count=\counter] \nodename in \namelist {
            \global\edef\numberofnodes{\counter}
            \node at (\nodename) [draw=none,name=hullnode\counter] {};
        }
        \node at (hullnode\numberofnodes) [name=hullnode0,draw=none] {};
        \pgfmathtruncatemacro\lastnumber{\numberofnodes+1}
        \node at (hullnode1) [name=hullnode\lastnumber,draw=none] {};
    },
    create hullnodes
]
($(hullnode1)!#2!-90:(hullnode0)$)
\foreach [
    evaluate=\currentnode as \previousnode using \currentnode-1,
    evaluate=\currentnode as \nextnode using \currentnode+1
    ] \currentnode in {1,...,\numberofnodes} {
-- ($(hullnode\currentnode)!#2!-90:(hullnode\previousnode)$)
  let \p1 = ($(hullnode\currentnode)!#2!-90:(hullnode\previousnode) - (hullnode\currentnode)$),
    \n1 = {atan2(\y1,\x1)},  
    \p2 = ($(hullnode\currentnode)!#2!90:(hullnode\nextnode) - (hullnode\currentnode)$),
    \n2 = {atan2(\y2,\x2)},  
    \n{delta} = {-Mod(\n1-\n2,360)}
  in 
    {arc [start angle=\n1, delta angle=\n{delta}, radius=#2]}
}
-- cycle
} 
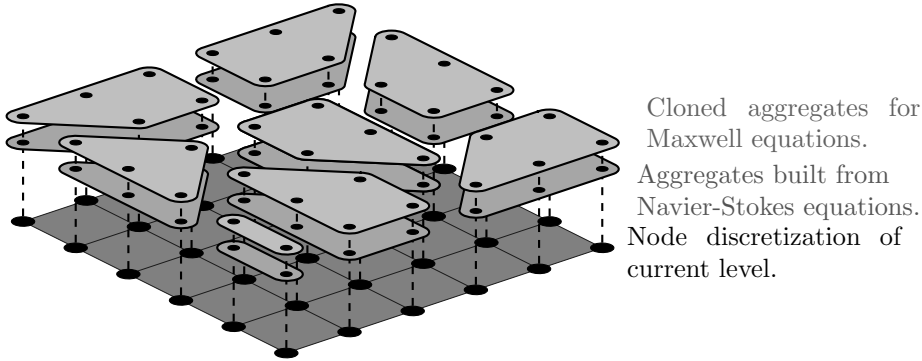
\begin{figure}[htbp]
    \centering
    \caption{Cloned aggregation strategy for volume-coupled multiphysics problems.}
    \label{fig:mgmphaggs}
\begin{tikzpicture}[scale=0.7]
\def\mycolorOne{\MphysicsA}  
\def\mycolorTwo{\MphysicsB}  
\def\mycolorThree{black}  
\def\mycolorFour{gray!60!white}  

\begin{scope}[every node/.append style={yslant=-0.5,xslant=1.2},yslant=-0.5,xslant=1.2]
\fill[color=gray] (0,0) rectangle (5,5);
\draw[step=1cm, line width=0.1mm, black!90!white] (0,0) grid (5,5);
  \foreach \x in {0,...,5}{
  \foreach \y in {0,...,5}{
        \node[draw,circle,inner sep=2pt,fill] at (\x,\y) (LOW\x\y) {};
            }
      }    
\end{scope}

\begin{scope}[shift={(0,1.5)},every node/.append style={yslant=-0.5,xslant=1.2},yslant=-0.5,xslant=1.2]
  \foreach \x in {0,...,5}{
  \foreach \y in {0,...,5}{
        \node[draw,circle,inner sep=2pt,fill] at (\x,\y) (UP\x\y) {};
            }
      } 
\end{scope}

\foreach \x in {0,...,5}{
  \foreach \y in {0,...,5}{
        \draw[thick,black,dashed] (LOW\x\y) -- (UP\x\y);
            }
      } 

\begin{scope}[shift={(0,1.5)},every node/.append style={yslant=-0.5,xslant=1.2},yslant=-0.5,xslant=1.2]       

   \draw[thick,black,fill=\mycolorOne,fill opacity=1.0] \convexpath{UP05,UP14,UP13,UP03}{2mm};
   \draw[thick,black,fill=\mycolorOne,fill opacity=1.0] \convexpath{UP15,UP35,UP34,UP24}{2mm};         
   \draw[thick,black,fill=\mycolorOne,fill opacity=1.0] \convexpath{UP45,UP55,UP53,UP44}{2mm};  
   \draw[thick,black,fill=\mycolorOne,fill opacity=1.0] \convexpath{UP22,UP23,UP43,UP32}{2mm};  
   \draw[thick,black,fill=\mycolorOne,fill opacity=1.0] \convexpath{UP22,UP23,UP43,UP32}{2mm};  
   \draw[thick,black,fill=\mycolorOne,fill opacity=1.0] \convexpath{UP00,UP02,UP12,UP11}{2mm}; 
   \draw[thick,black,fill=\mycolorOne,fill opacity=1.0] \convexpath{UP10,UP21,UP30}{2mm}; 
   \draw[thick,black,fill=\mycolorOne,fill opacity=1.0] \convexpath{UP42,UP52,UP51,UP31}{2mm};       
   \draw[thick,black,fill=\mycolorOne,fill opacity=1.0] \convexpath{UP40,UP50}{2mm};             

  \foreach \x in {0,...,5}{
  \foreach \y in {0,...,5}{
        \node[draw,circle,inner sep=1pt,fill] at (\x,\y) {};
            }
      } 
      
\end{scope}

\begin{scope}[shift={(0,2.0)},every node/.append style={yslant=-0.5,xslant=1.2},yslant=-0.5,xslant=1.2]
  \foreach \x in {0,...,5}{
  \foreach \y in {0,...,5}{
        \node[draw,circle,inner sep=2pt,fill] at (\x,\y) (UPB\x\y) {};
            }
      } 
\end{scope}

\foreach \x in {0,...,5}{
  \foreach \y in {0,...,5}{
        \draw[thick,black,dashed] (UP\x\y) -- (UPB\x\y);
            }
      } 

\begin{scope}[shift={(0,2.0)},every node/.append style={yslant=-0.5,xslant=1.2},yslant=-0.5,xslant=1.2]       

   \draw[thick,black,fill=\mycolorTwo,fill opacity=1.0] \convexpath{UPB05,UPB14,UPB13,UPB03}{2mm};
   \draw[thick,black,fill=\mycolorTwo,fill opacity=1.0] \convexpath{UPB15,UPB35,UPB34,UPB24}{2mm};         
   \draw[thick,black,fill=\mycolorTwo,fill opacity=1.0] \convexpath{UPB45,UPB55,UPB53,UPB44}{2mm};  
   \draw[thick,black,fill=\mycolorTwo,fill opacity=1.0] \convexpath{UPB22,UPB23,UPB43,UPB32}{2mm};  
   \draw[thick,black,fill=\mycolorTwo,fill opacity=1.0] \convexpath{UPB22,UPB23,UPB43,UPB32}{2mm};  
   \draw[thick,black,fill=\mycolorTwo,fill opacity=1.0] \convexpath{UPB00,UPB02,UPB12,UPB11}{2mm}; 
   \draw[thick,black,fill=\mycolorTwo,fill opacity=1.0] \convexpath{UPB10,UPB21,UPB30}{2mm}; 
   \draw[thick,black,fill=\mycolorTwo,fill opacity=1.0] \convexpath{UPB42,UPB52,UPB51,UPB31}{2mm};       
   \draw[thick,black,fill=\mycolorTwo,fill opacity=1.0] \convexpath{UPB40,UPB50}{2mm};             

  \foreach \x in {0,...,5}{
  \foreach \y in {0,...,5}{
        \node[draw,circle,inner sep=1pt,fill] at (\x,\y) {};
            }
      } 
      
\end{scope}

\node[below right= 0.1cm and 0.35cm of UPB55] {\color{\mycolorOne!50!black}{\begin{minipage}{0.29\textwidth}Aggregates built from \newline Navier-Stokes equations.\end{minipage}}};
\node[above right= 0.1cm and 0.3cm of UP55] {\color{\mycolorTwo!50!black}{\begin{minipage}{0.28\textwidth}Cloned aggregates for Maxwell equations.\end{minipage}}};

\node[right= 0.1cm of LOW55] {\color{black}{\begin{minipage}{0.28\textwidth}Node discretization of current level.\end{minipage}}};
\end{tikzpicture}%
\end{figure}
As discussed later, this type of multigrid adaptation is relatively
straight forward with {\sf MueLu} using an {\sf XML} input file.
In addition to consistently coarsening both fluid and electromagnetic variables,
the reuse or cloning of aggregates saves a modest amount of time within the 
multigrid setup phase. Additional reductions in the setup cost can be achieved if one
defines $\Restr_{11} = \Restr_{00}$ and $\Prol_{11} = \Prol_{00}$. Again, this
is relatively easy to do within the {\sf MueLu} framework, though not necessary.
In our experiments, piecewise constant interpolation is the basis for all grid
transfers. This simple grid transfer choice is more robust for highly convective flows.

\subsubsection{Block smoother for the MHD system}


A block Gauss-Seidel (BGS) iteration forms the basis of the multigrid smoother using the 
notion of blocks already introduced via the $2\times 2$ decomposition. A 
standard block Gauss-Seidel iteration would {\it solve} for an entire block of 
unknowns simultaneously, recompute residuals, and then {\it solve} for the
other block of equations simultaneously. This corresponds to alternating
between a Navier-Stokes sub-block solve and a Maxwell sub-block solve while 
performing residual updates after each sub-solve. Effectively, the 
 $2\times 2$ block perspective emphasizes the coupling within the Navier-Stokes
block and within the Maxwell part as each of these blocks correspond to a saddle-point
system due to the presence of constraint equations. That is, the difficulties
associated with saddle-point systems will be addressed by an approximate 
sub-block solve. Further, the block Gauss-Seidel iteration considers 
the Navier-Stokes block and the Maxwell block equally important as the 
iteration simply alternates equally between the two sub-solves. As the exact
solution of each sub-system is computationally expensive, domain decomposition with ILU(0) is used
to generate approximate sub-solutions associating one domain with each computing core. 
While not necessarily an optimal smoother, we 
have found that ILU(0) is often effective for the saddle-point systems associated
with incompressible flow.
While ILU(0) could be applied to the 
entire $2\times 2$ matrix, this involves significantly more computation and memory during
the setup and apply phases of the smoother. This more expensive process is generally 
not needed as the coupling between the block sub-systems is often much
less important than the coupling within the sub-blocks, though there could 
be MHD situations with significant cross-coupling between Navier-Stokes and Maxwell
that would warrant an ILU smoother applied to the whole system.

\begin{algorithm}
\caption{Damped blocked Gauss-Seidel smoother}
\label{alg:bgsiteration}
\begin{algorithmic}[1]
\REQUIRE{$\A,\rhsv,\bgsdamping,\#\textnormal{sweeps}$}
\STATE{Set initial guess: $\solv:=\textbf{0}$}
\FOR{$s=0$, $s<\#\textnormal{sweeps}$}
\STATE{\textit{\% Calculate update for Navier-Stokes part}}
\STATE{\label{algResA}Calculate residual: $\resvO:=
\rhsvO-\AOO \solvO - \AOl  \solvl$}
\STATE{\label{algSolveA}Solve approximately $\AOO\widetilde{\solvO}=\resvO$ for $\widetilde\solvO$}
\STATE{\label{algScalingA}Update intermediate solution: $\solvO:= \solvO + \bgsdamping~ \widetilde\solvO$}
\STATE{\textit{\% Calculate update for Maxwell part}}
\STATE{\label{algResB}Calculate residual: $\resvl:=\rhsvl- \AlO \solvO - \All \solvl$}
\STATE{\label{algSolveB}Solve approximately $\All\widetilde{\solvl}=\resvl$ for $\widetilde\solvl$}
\STATE{\label{algScalingB}Update intermediate solution: $\solvl:=\solvl + \bgsdamping~ \widetilde\solvl$} 
\ENDFOR
\RETURN{$\solv:=\begin{bmatrix}\solvO\\\solvl\end{bmatrix}$}
\end{algorithmic}
\end{algorithm}

Algorithm \ref{alg:bgsiteration} shows the outline of the damped Gauss-Seidel block smoothing algorithm. First, an approximate solution update $\widetilde\solvO$ of the Navier-Stokes part is built in line \ref{algSolveA} of Algorithm \ref{alg:bgsiteration} and then scaled by a damping parameter $\bgsdamping$ in line \ref{algScalingA}. Similarly, an approximate solution update $\widetilde\solvl$ is built for the Maxwell part and scaled by the same damping parameter $\bgsdamping$ in lines \ref{algSolveB} and \ref{algScalingB}. Please note, that the residual calculation in line \ref{algResB} employs the intermediate solution update from line \ref{algScalingA}. Similarly, an intermediate solution is used for the residual calculation in line \ref{algResA}, if we apply more than one sweep with the block Gauss-Seidel smoothing algorithm.

As one can see from Algorithm \ref{alg:bgsiteration} we only need approximate inverses of the diagonal blocks $\AOO$ and $\All$. A flexible implementation allows one to choose appropriate local smoothing methods to approximately invert the blocks $\AOO$ and $\All$. The coupling is guaranteed by the off-diagonal blocks $\AOl$ and $\AlO$ in the residual calculations in lines \ref{algResA} and \ref{algResB}.

\subsection{Software}

Even though this paper demonstrates the multigrid approach for certain MHD formulations, it is clear that the concept is more general. Due to the problem-dependent 
nature of preconditioning multiphysics problems, 
it is logical to provide a general software framework that helps design application-specific preconditioning methods.
The proposed methods in this paper are implemented using the next-generation multigrid framework MueLu from the Trilinos software libraries.
In contrast to other publications like \cite{verdugo2016} 
which share the same core idea of a general software framework, MueLu is publicly available through the Trilinos library. Furthermore, it is fully embedded in the Trilinos software stack and has full native access to all features provided by the other Trilinos packages. It aims at next-generation HPC platforms and automatically benefits from all performance improvements in the underlying core linear algebra packages.

The core design concept of MueLu is based on building blocks that can be combined to construct complex preconditioner layouts. Figure \ref{fig:standardsetup} 
shows an algorithmic layout to build standard aggregation-based transfer operators $\Prol$ and $\Restr$, the level smoother $\mathcal{S}$ and a new 
coarse level operator $\A_{\ell+1}$. Further coarse levels follow by recursively applying a multigrid setup procedure to the  coarse level operators.
Each building block, denoted by a rectangular block in Figure \ref{fig:standardsetup}, processes certain input data and produces output data which serves as input for the downstream building blocks. For application-specific adaptions it is usually sufficient to replace specific algorithmic components while the majority of building blocks can be reused.

Figure \ref{fig:standardsetup} represents the process of defining coarse multigrid levels during the setup phase for a single-field multigrid method.
%
%
Figure \ref{fig:multiphysicssetup} shows the algorithmic design of the setup phase to create a $2\times 2$ block transfer operator for a volume-coupled problem as 
discussed earlier. 
Here, aggregates built from block $\AOO$ are re-used instead of re-building them from build block $\All$

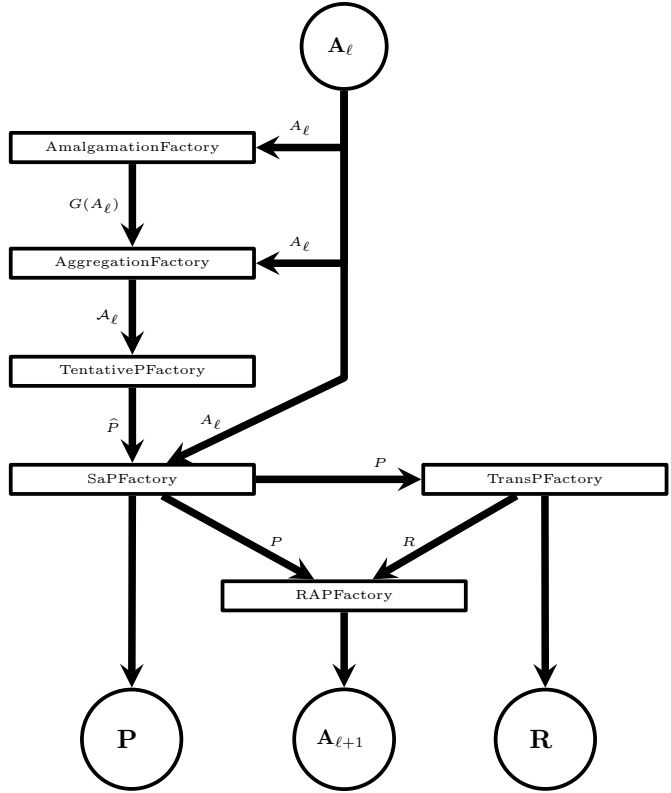
\begin{figure}
\caption{Algorithmic layout for standard preconditioner setup}
\label{fig:standardsetup}
\centering
	\begin{tikzpicture}[line join=round]
			\begin{scope}[shift={(0,0)}]
			\tikzstyle{style_var}=[anchor=center,circle, draw=black, fill=white, text centered, anchor=north, text=black, text width=0.8cm, line width=0.5mm]
			\tikzstyle{style_var_output}=[anchor=center,circle, draw=black, fill=white, text centered, anchor=north, text=black, text width=1.0cm, line width=0.5mm]

			\tikzstyle{myarrow}=[->, >=stealth, line width=1mm]

		\begin{scope}[shift={(0,3)}]
		\tikzstyle{style_factory}=[anchor=center,rectangle, draw=black, fill=white, text centered, anchor=north, text=black, text width=3cm, line width=0.5mm]
		\node (Afffine) [style_var]  { \scalebox{0.8}{$\A_{\ell}$ }};
		\node[node distance=6.5] (RAPFactory) [style_factory, below=of Afffine]
		{
		  \tiny RAPFactory
		};
		\node (Acoarse) [style_var_output, below=of RAPFactory]  { \scalebox{0.8}{$\A_{\ell+1}$ }};
		\node (Pcoarse) [node distance=1.45,style_var_output, left=of Acoarse]  { \scalebox{1.0}{$\Prol$ }};
		\node (Rcoarse) [node distance=1.3,style_var_output, right=of Acoarse]  { \scalebox{1.0}{$\Restr$ }};	
		\end{scope}

		\tikzstyle{style_factory}=[anchor=center,rectangle, draw=black, fill=white, text centered, anchor=north, text=black, text width=3cm, line width=0.5mm]
		\tikzstyle{style_blockedfactory}=[anchor=center,rectangle, draw=black, fill=white, text centered, anchor=north, text=black, text width=3cm, line width=0.5mm]

		\node (lHelper) [node distance=2.2, inner sep=0pt, left=of Afffine] {};		
		\node[node distance=1.1] (AmalgamationFactory) [style_factory, below=of lHelper]
		{
		  \tiny AmalgamationFactory
		};
		\node[node distance=1.1] (AggregationFactory) [style_factory, below=of AmalgamationFactory]
		{
		  \tiny AggregationFactory
		};
		\node[node distance=1.0] (Ptent) [style_factory, below=of AggregationFactory]
		{
		  \tiny TentativePFactory
		};
		\coordinate (Ahelper) at (0,-2);
		\node[node distance=1.0] (PSmoothing) [style_factory, below=of Ptent]
		{
		  \tiny SaPFactory
		};
		\node[node distance=2.2] (Restrictor) [style_factory, right=of PSmoothing]
		{
		  \tiny TransPFactory
		};

		\draw[myarrow,->] (Afffine) -- (Ahelper) -- (PSmoothing) node [above, near end]{\tiny $A_\ell$};
		\draw[myarrow,->] (Afffine) |- (AggregationFactory) node [above, near end]{\tiny $A_\ell$};
		\draw[myarrow,->] (Afffine) |- (AmalgamationFactory) node [above, near end]{\tiny $A_\ell$};
		\draw[myarrow,->] (AmalgamationFactory) -- (AggregationFactory) node [left, midway] {\tiny $G(A_\ell)$};
		\draw[myarrow,->] (AggregationFactory) -- (Ptent) node [left, midway] {\tiny $\mathcal{A}_\ell$};
		\draw[myarrow,->] (Ptent) -- (PSmoothing) node [left, midway] {\tiny $\widehat{P}$};
		\draw[myarrow,->] (PSmoothing) -- node [above,near end] {\tiny $P$} (Restrictor);
	
		\draw[myarrow,->] (Restrictor) -- (Rcoarse);	
		\draw[myarrow,->] (PSmoothing) -- (Pcoarse);			
		\draw[myarrow,->] (RAPFactory) -- (Acoarse);			
		\draw[myarrow,->] (PSmoothing) -- node [above,near end] {\tiny $P$} (RAPFactory);			
		\draw[myarrow,->] (Restrictor) -- node [above,near end] {\tiny $R$} (RAPFactory);					
		\end{scope} 
		
	\end{tikzpicture}%
\end{figure}
	
\begin{figure}
\caption{Algorithmic multigrid setup layout for $2\times 2$ block multigrid method}
\label{fig:multiphysicssetup}
\centering
    \scalebox{0.66}{%
	\begin{tikzpicture}
		\tikzstyle{style_var}=[anchor=center,circle, draw=black, fill=white, text centered, anchor=north, text=black, text width=0.8cm, line width=0.5mm]
		\tikzstyle{style_var_output}=[anchor=center,circle, draw=black, fill=white, text centered, anchor=north, text=black, text width=0.8cm, line width=0.5mm]

		\tikzstyle{myarrow}=[->, >=stealth, line width=1mm]
		\tikzstyle{style_factory}=[anchor=center,rectangle, draw=black, fill=white, text centered, anchor=north, text=black, text width=3cm, line width=0.5mm]

		\node (Afffine) [style_var] at (0,2) { \scalebox{1.0}{$\A_{\ell}$ }};
		\node (RAPFact) [style_factory] at (0,-8) { \tiny RAPFactory};

		\begin{scope}[shift={(-6,0)}]%
		\tikzstyle{myarrow}=[->, >=stealth, line width=1mm]
		\tikzstyle{style_factory}=[anchor=center,rectangle, draw=black, fill=\MphysicsA, text centered, anchor=north, text=black, text width=3cm, line width=0.5mm]
		\tikzstyle{style_blockedfactory}=[anchor=center,rectangle, draw=black, fill=white, text centered, anchor=north, text=black, text width=3cm, line width=0.5mm]
		\node (Afine) [style_factory]  { \tiny SubBlockAFactory(0,0)};
		\node (Pfine) [node distance=1,right=of Afine,opacity=0.7,text=gray]     {};
		\draw[myarrow,->] (Afffine) -| (Afine);		
		\node[node distance=1.1] (AmalgamationFactory) [style_factory, below=of Pfine]
		{
		  \tiny AmalgamationFactory
		};
		\node[node distance=1.1] (AggregationFactory) [style_factory, below=of AmalgamationFactory]
		{
		  \tiny AggregationFactory
		};
		\node[node distance=1.0] (Ptent) [style_factory, below=of AggregationFactory]
		{
		  \tiny TentativePFactory
		};
		\node[node distance=1] (RestrictorA) [style_factory, left=of Ptent]
		{
		  \tiny TransPFactory
		};
		\draw[myarrow,->] (Afine) |- (AggregationFactory) node [above, near end]{\tiny $A_\ell$};
		\draw[myarrow,->] (Afine) |- (AmalgamationFactory) node [above, near end]{\tiny $A_\ell$};
		\draw[myarrow,->] (AmalgamationFactory) -- (AggregationFactory) node [left, midway] {\tiny $G(A_\ell)$};
		\draw[myarrow,->] (AggregationFactory) -- (Ptent) node [left, midway] {\tiny $\mathcal{A}_\ell$};
		\draw[myarrow,->] (Ptent) -- node [above,near end] {\tiny $P$} (RestrictorA);
		\node (BlockedProlongator) [style_blockedfactory,node distance=1.5, below=of Ptent]
		{
		  \tiny BlockedTransferFactory
		};
		\draw[myarrow,->] (Ptent) -- node [left,near end] {\tiny $P$} (BlockedProlongator);		
		\end{scope}
		
		\begin{scope}[shift={(6,0)}]%
		\tikzstyle{myarrow}=[->, >=stealth, line width=1mm]
		\tikzstyle{style_factory}=[anchor=center,rectangle, draw=black, fill=\MphysicsB, text centered, anchor=north, text=black, text width=3cm, line width=0.5mm]
		\tikzstyle{style_blockedfactory}=[anchor=center,rectangle, draw=black, fill=white, text centered, anchor=north, text=black, text width=3cm, line width=0.5mm]
		\node (Afine) [style_factory]  { \tiny SubBlockAFactory(1,1)};
		\node (Pfine) [node distance=1,left=of Afine,opacity=0.7,text=gray]     {};
		\draw[myarrow,->] (Afffine) -| (Afine);		
		\node[node distance=1.1] (AmalgamationFactory) [style_factory, below=of Pfine]
		{
		  \tiny AmalgamationFactory
		};
		\node[node distance=2.6] (PtentB) [style_factory, below=of AmalgamationFactory]
		{
		  \tiny TentativePFactory
		};
		\node[node distance=1] (RestrictorB) [style_factory, right=of PtentB]
		{
		  \tiny TransPFactory
		};
		\draw[myarrow,->] (Afine) |- (AmalgamationFactory) node [above, near end]{\tiny $A_\ell$};
		\draw[myarrow,->] (AmalgamationFactory) -- (PtentB) node [right, midway] {\tiny $G(A_\ell)$};
		\draw[myarrow,->] (PtentB) -- node [above,near end] {\tiny $P$} (RestrictorB);
		\draw[myarrow,->] (AggregationFactory) -| (PtentB) node [left, near end] {\tiny $\mathcal{A}_\ell$};
		\draw[myarrow,->,color=red] (AggregationFactory) -| (PtentB) node [left, near end] {\tiny $\mathcal{A}_\ell$};
		\node (BlockedRestrictor) [style_blockedfactory,node distance=1.5, below=of PtentB]
		{
		  \tiny BlockedTransferFactory
		};
		\draw[myarrow,->] (PtentB) -- node [below,near end] {\tiny $P$} (BlockedProlongator);		
		\draw[myarrow,->] (RestrictorB) -- node [below] {\tiny $R$} (BlockedRestrictor);			
		\draw[myarrow,->] (RestrictorA) -- node [below] {\tiny $R$} (BlockedRestrictor);	
		\end{scope}
		\draw[myarrow,->] (BlockedRestrictor) -- node [above,near end] {\tiny $R$} (RAPFact);			
		\draw[myarrow,->] (BlockedProlongator) -- node [above,near end] {\tiny $P$} (RAPFact);	
		\draw[myarrow,->] (Afffine) -- node [left] {\tiny $\A_\ell$} (RAPFact);	

		\node (Acoarse) [style_var,below=of RAPFact,node distance=1] { \scalebox{0.8}{$\A_{\ell+1}$ }};
		\node (P) [style_var,node distance=2.6,below=of BlockedProlongator] { \scalebox{1.0}{$\Prol$ }};		
		\node (R) [style_var,node distance=2.6,below=of BlockedRestrictor] { \scalebox{1.0}{$\Restr$ }};		
		
		\draw[myarrow,->] (RAPFact) -- (Acoarse);			
		\draw[myarrow,->] (BlockedProlongator) -- (P);	
		\draw[myarrow,->] (BlockedRestrictor) -- (R);	

        \end{tikzpicture}%
        }
\end{figure}
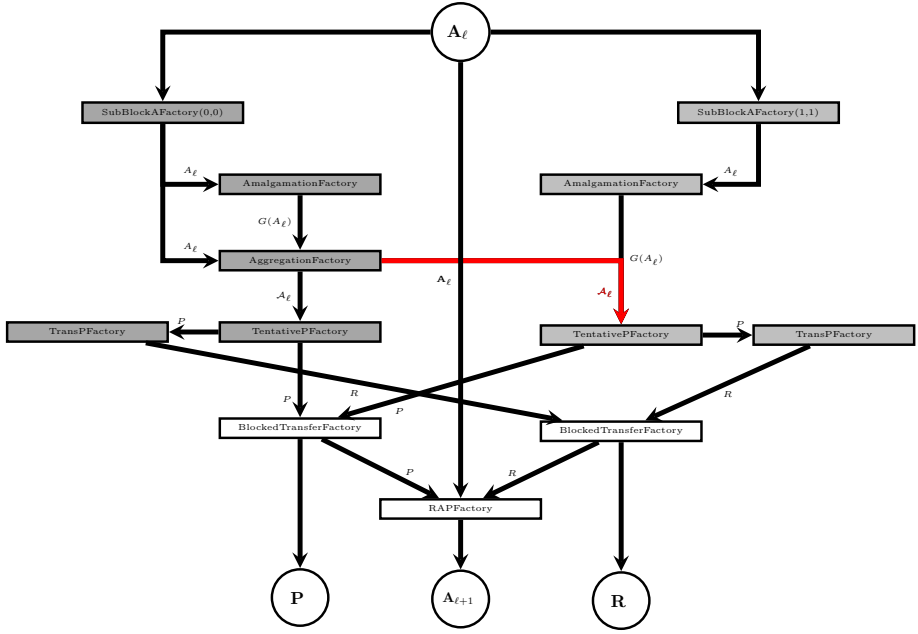

\section{Experimental results}
\label{sec:experiments}

%
%

%
%
%
%
%

%
%
%

\providecommand{\tablelegendnew}{%
\renewcommand{\arraystretch}{1.0}%
\medskip{%
\begin{center}%
\footnotesize%
Legend:%
\begin{tabular}{cl}%
$n_N$ & Accumulated number of all nonlinear iterations \\%
$n_L$ & Accumulated number of all linear iterations \\%
$t_{\Sigma}$ & Solver time (setup + iteration phase) \\%
\end{tabular}%
\end{center}%
}%
}

\subsection{MHD generator}
\label{sec:mhdgen}
This problem is a steady-state MHD duct flow configuration representing an idealized MHD generator  
where an electrical current is induced by pumping a conducting fluid (mechanical work) through an externally applied vertical magnetic field \cite{shadid2016scalable}. 
The bending of the magnetic field lines produces a horizontal electrical current.
The geometric domain for this problem is a square cross-sectional duct of dimensions $[0,15]\times[0,1]\times[0,1]$.
The velocity boundary conditions are set with Dirichlet inlet velocity of $\vel = (\vels,0,0)$, no-slip on the top, bottom and sides of the channel, and natural boundary conditions on the outflow. 
The magnetic field on the top and bottom boundaries is specified as  $\magn=(0,\magns_y^{\textnormal{gen}},0)$ where 
\begin{equation*}
\magns_y^{\textnormal{gen}}   =  \frac{1}{2} \magns_0 \biggl[ \tanh\Bigl(\frac{x - x_{\textnormal{on}}}{\delta}\Bigr) -  \tanh\Bigl(\frac{x - x_{\textnormal{off}}}{\delta}\Bigr) \biggr].
\end{equation*} 
Here, $B_0$ is the strength of the field and $\delta$ is a measure of the transition length-scale for application of the field. 
The inlet, outlet and sides  are perfect conductors with $\magn \cdot \hat {\bf n} = {\bf 0}$ and ${\bf E} \times \hat {\bf n} = {\bf 0}$, where $\hat {\bf n}$  is the outward facing unit normal vector. 
Zero Dirichlet boundary conditions are applied on all surfaces for the Lagrange multiplier.
The problem is defined by three non-dimensional parameters: 
the Reynolds number $Re = \rho u L / \mu$,
the magnetic Reynolds number $ Re_m = \mu_0 u L/ \eta$, 
and the Hartmann number $Ha = B_0 L / \sqrt{ \rho \nu \eta}$.
Here, $u$ is the maximum x-direction velocity.
The parameters in this problem are taken to be $\vels = 1.0$, $\dens =1$, $B_0 = 3.354$,  $\magnperm =1$, 
$\resistivity = 1$, $x_{\textnormal{on}} = 4.0$, $x_{\textnormal{off}} = 6.0$, and $\delta = 0.5$. 
The linear solver used was non-restarted GMRES to reach a relative tolerance of $10^{-3}$.
As a reference preconditioner we use a fully-coupled or non-blocked AMG method (FC-AMG),
where the relaxation method, Additive Schwarz with overlap one,
is applied to the entire system that includes the off-diagonal coupling between the fluids and magnetics (see~\cite{lin2010AMGjcp}).
For the blocked variant presented in this manuscript (Section~\ref{sec:blockAMG}), we use BGS as a relaxation method (AMG(BGS).
The approximate solves for the sub-blocks are handled with Additive Schwarz with overlap one.
The coarse grid for both the FC-AMG and AMG(BGS) is solved directly.

In the first study, we investigate the number of iterations and solution time as a function of 
 the block Gauss-Seidel (BGS) damping parameter $\bgsdamping$ for various viscosities\\ 
$\visc \in \{0.006,0.007,0.008,0.01,0.02,0.04 \}$, which effectively sets the range for the non-dimensional parameters: 
$25 \le Re \le 167$, $17 \le Ha \le 43$, and $Re_m = 1$.
\begin{figure}[htbp]
\centering
\caption{Solver performance for AMG(BGS) preconditioner with 1 BGS coupling iteration versus the fully coupled AMG(ILU) preconditioner for MHD generator example on a $240\times 16 \times 16$ mesh using $32$ processors.}
\label{fig:mhdgenmode1}
\begin{subfigure}{0.45\textwidth}
\includegraphics[width=0.98\textwidth]{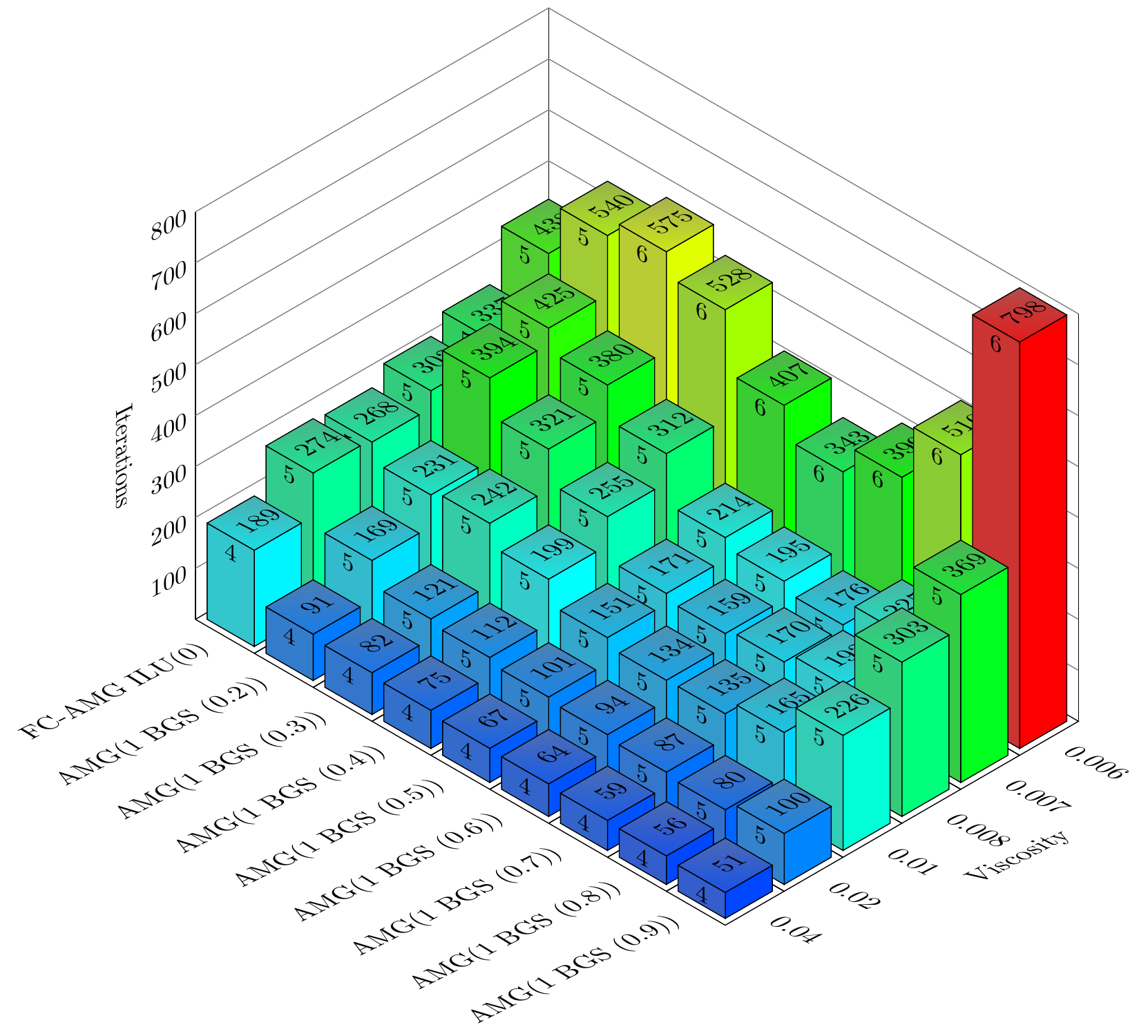}
\caption{Accumulated number of linear iterations depending on BGS damping parameter and viscosity. The numbers on the sides of the column denote the number of nonlinear iterations to solve the system.}
\label{fig:mhdgenmode1a}
\end{subfigure}~~%
\begin{subfigure}{0.45\textwidth}
\includegraphics[width=0.98\textwidth]{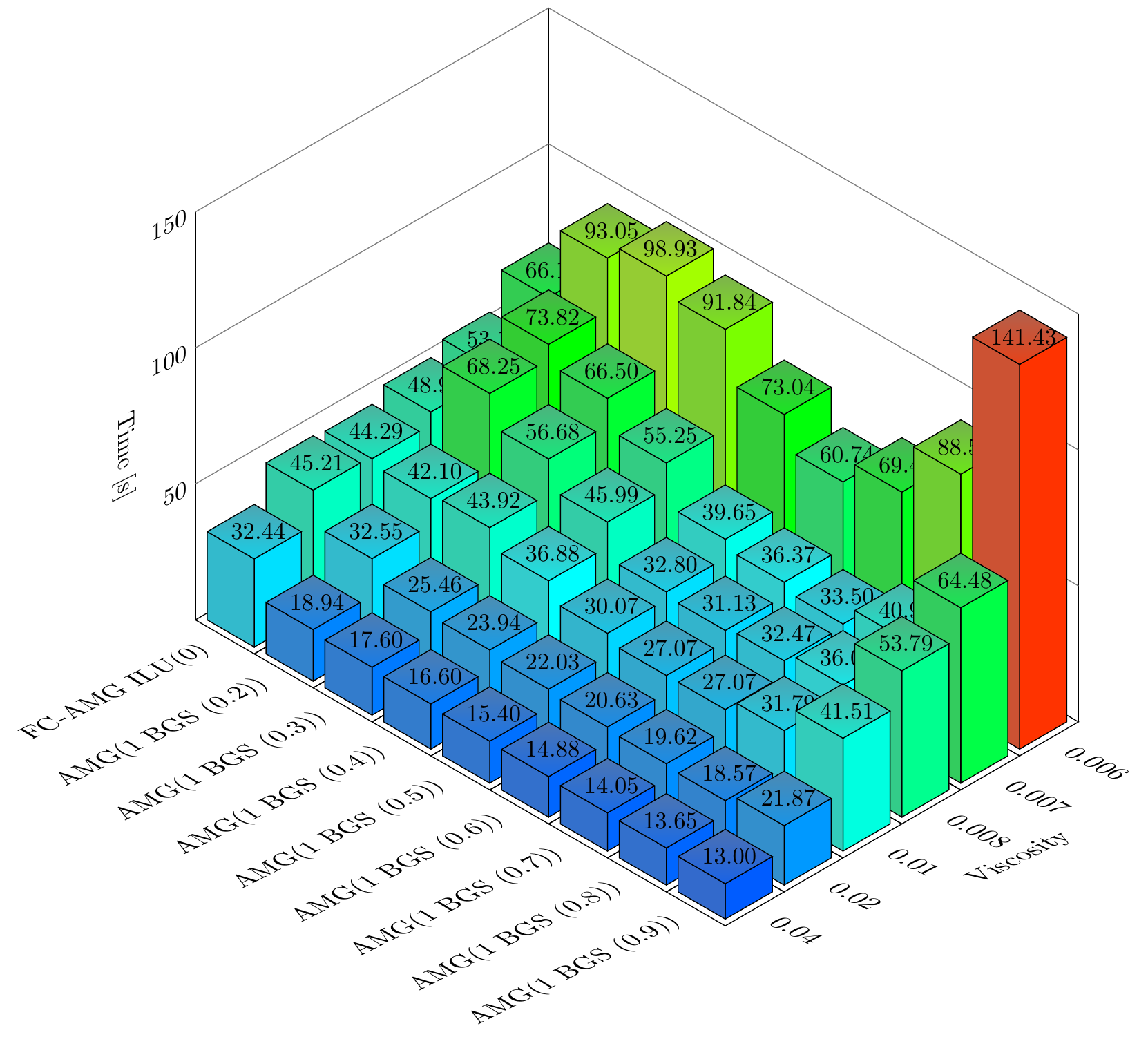}
\caption{Solver timings (setup + iteration phase) depending on BGS damping parameter and viscosity.}
\label{fig:mhdgenmode1b}
\end{subfigure}
\end{figure}

Figure \ref{fig:mhdgenmode1} shows the accumulated number of linear iterations and solver timings for the static MHD generator problem on a $240\times 16\times 16$ mesh using $32$ processors (Intel Broadwell E5-2695 (2.1 GHz), 1 cluster node with $128$ GB RAM). 
The numbers on the side of the columns in Figure \ref{fig:mhdgenmode1a} denote the number of nonlinear iterations. The number on the z-axis 
represents the accumulated number of linear iterations.
Figure \ref{fig:mhdgenmode1b} displays the accumulated solver timings (setup and iteration phase) of the corresponding preconditioning variants.
The timings are averaged over 5 simulations.
As one can see from Figure \ref{fig:mhdgenmode1}, the right choice of the BGS damping parameter is crucial for smaller viscosities. Even though the optimal damping parameter depends on the problem, in practice a damping parameter $\bgsdamping$ near $0.5$ seems to work well.
Figure \ref{fig:mhdgenmode2} shows the accumulated linear iterations and the solver times for the MHD generator example on a finer $480\times 32\times 32$ mesh using $256$ processors
(Intel Broadwell E5-2695 (2.1 GHz), 8 cluster nodes with 128 GB RAM each, Intel Omni-Path high speed interconnect).
One can see in Figure \ref{fig:mhdgenmode2a} that a higher number of BGS iterations reduces the overall number of linear iterations,
but not enough to compensate for the higher costs per iteration (see Figure \ref{fig:mhdgenmode2b}).
For reference, the fully coupled approach 
is also shown in the back row (labeled as \textit{FC-AMG ILU(0)}).
This option does not require the multiphysics framework, but it generally takes more
time and iterations than the multiphysics solver, varying somewhere between
$2x$ and $3x$ slower. 

\begin{figure}[htbp]
\centering
\caption{Solver performance for AMG(BGS) preconditioner with fixed damping parameter and 1 or 2 BGS coupling iterations versus the fully coupled AMG(ILU) preconditioner for MHD generator example on a $480\times 32 \times 32$ mesh using $256$ processors.}
\label{fig:mhdgenmode2}
\begin{subfigure}{0.45\textwidth}
\includegraphics[width=0.98\textwidth]{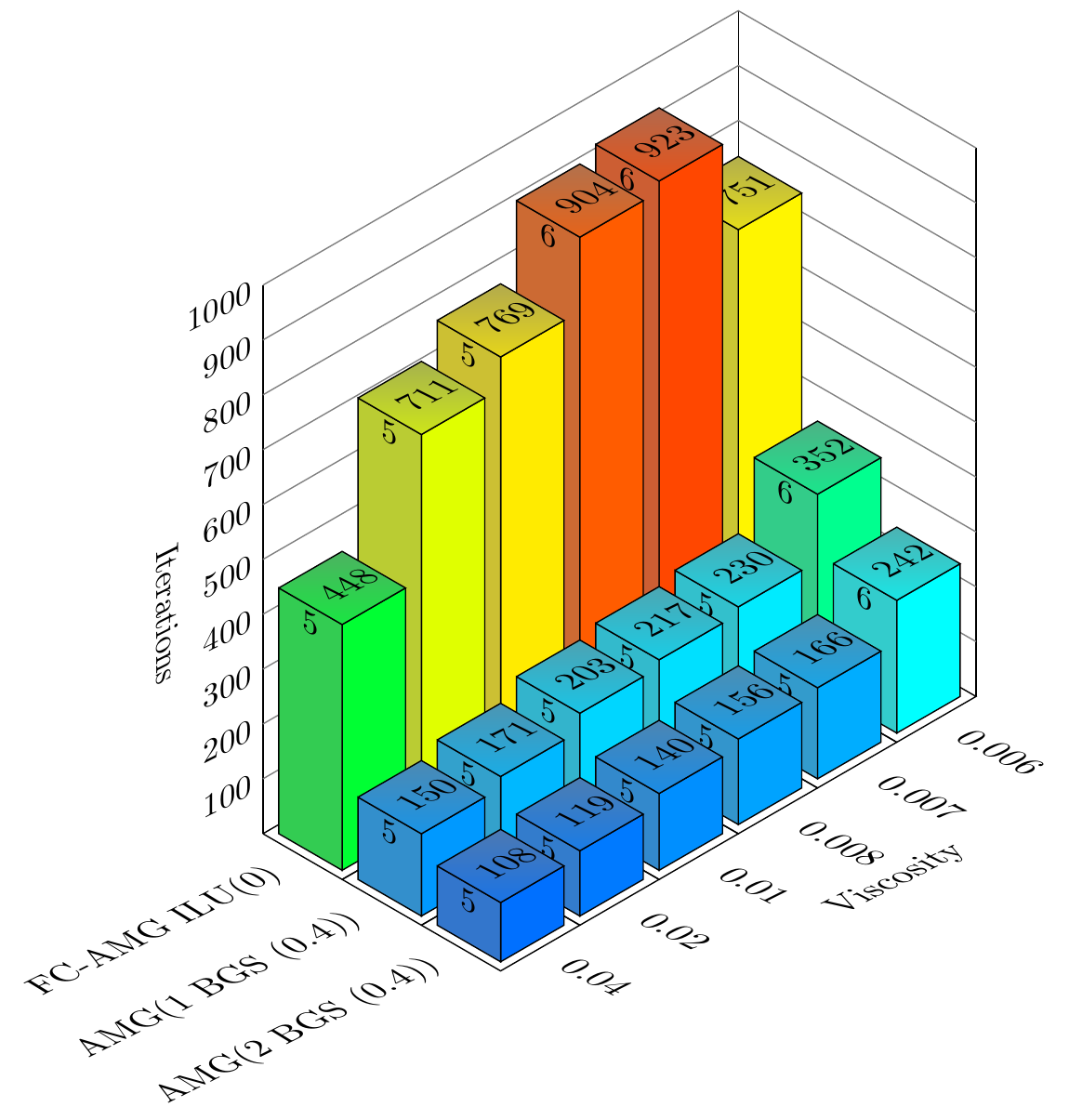}
\caption{Accumulated number of linear iterations depending on BGS iterations and viscosity. The numbers on the sides of the column denote the number of nonlinear iterations to solve the system.}
\label{fig:mhdgenmode2a}
\end{subfigure}~~%
\begin{subfigure}{0.45\textwidth}
\includegraphics[width=0.98\textwidth]{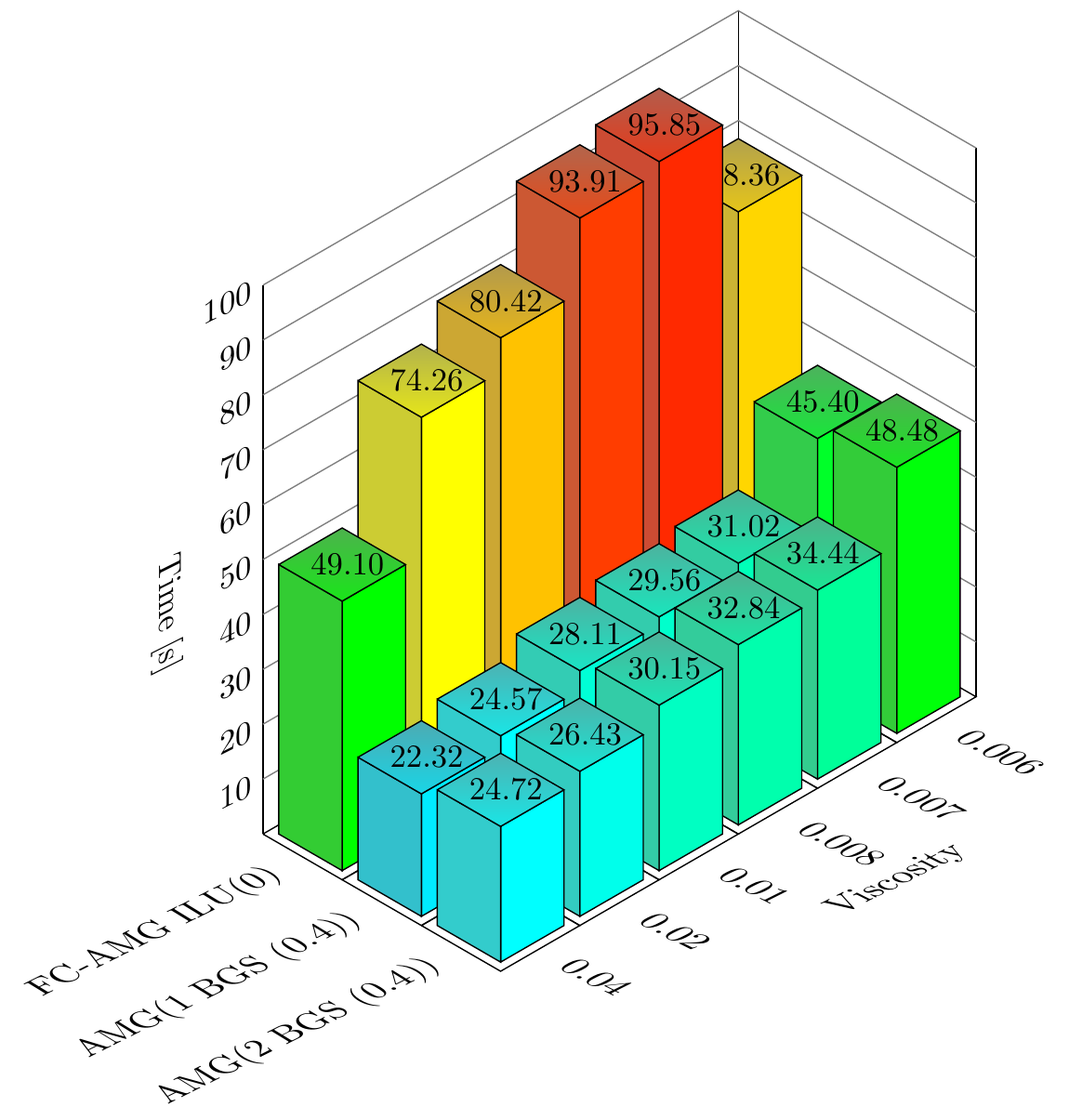}
\caption{Solver timings (setup + iteration phase) depending on BGS iterations parameter and viscosity.}
\label{fig:mhdgenmode2b}
\end{subfigure}
\end{figure}

It should be noted that the scaling of the iteration counts with respect to mesh refinement can be improved by increasing the linear solver tolerance for both the FC-AMG and AMG(BGS) methods. In experiments not shown here using the same solver settings but tightening the tolerance to $10^{-6}$, the FC-AMG total number of linear iterations increased by less than $10\%$ while the AMG(BGS) iterations decreased for the same two meshes considered in Figure \ref{fig:mhdgenmode1} and \ref{fig:mhdgenmode2}. While the scaling is better, the overall solution time with these tighter tolerances is longer and therefore it is not considered further in this paper since for engineering applications our primary goal is to minimize the solver time and not the iteration count.

\providecommand{\tablelegend}{%
\renewcommand{\arraystretch}{1.0}%
\medskip{%
\begin{center}%
\footnotesize%
Legend:%
\begin{tabular}{cl}%
$n_T$ & Number of time steps \\%
$n_N$ & Accumulated number of all nonlinear iterations \\%
$n_L$ & Accumulated number of all linear iterations \\%
$t_{Se}$ & Multigrid setup time \\%
$t_{So}$ & Multigrid solution time (iteration phase)\\%
$t_{\Sigma}$ & Solver time (setup + iteration phase) \\%
\end{tabular}%
\end{center}%
}%
}

\providecommand{\tablelegendmixed}{%
\renewcommand{\arraystretch}{1.0}%
\medskip{%
\begin{center}%
\footnotesize%
Legend:%
\begin{tabular}{cl}%
$n_N$ & Accumulated number of all nonlinear iterations \\%
$n_L$ & Accumulated number of all linear iterations \\%
$t_{Se}$ & Multigrid setup time \\%
$t_{So}$ & Multigrid solution time (iteration phase)\\%
\end{tabular}%
\end{center}%
}%
}

\subsection{Hydromagnetic Kevin-Helmholtz (HMKH)}

\label{sec:HMKH}
A hydromagnetic Kelvin-Helmholtz unstable shear layer problem is a configuration used to study magnetic reconnection \cite{shadid2016scalable} and
is posed in a domain of $[0,4]\times[-2,2]\times[0,2]$. It is described by an initial condition defined by two counter flowing conducting fluid streams with constant velocities $\vel (x,y>0,z,0) = (5,0,0)$ and $\vel (x,y<0,z,0) = (-5,0,0)$ and a Harris sheet sheared magnetic field configuration given by
\begin{equation}
\magn (x,y,z,0) = \left(0,\magns_0 \tanh(y/\delta),0\right).
\end{equation}
The boundary conditions are periodic on the right and left, as well as the front and back. 
The top and bottom are impenetrable for the fluid velocity, and the magnetic field is defined by the Harris sheet.
The magnetic Lagrange multiplier is taken as zero on all boundaries. 
The parameters in this problem are $\dens =1$, $\magnperm =1$, 
$\visc = 10^{-4}$,$\resistivity = 10^{-4}$, $\magns_0 = 0.3333$, $\delta = 0.1$ to produce $Re =  5 \times 10^4$, 
 $Re_m = 5 \times 10^4$, and an Alfven velocity, $u_A = B_0 / \sqrt{\rho \mu_0} = 0.333$, resulting in an Alfvenic Mach number $M_A = u /u_A = 15$. For these non-dimensional parameters, the shear layer is
Kelvin-Helmholtz unstable and forms a vortex sheet that evolves with time and undergoes 
thin current sheet formation, vortex rollup and merging. 
Figure \ref{fig:hmkh} shows the unstable shear layer  evolving from 
smaller vortices to a larger vortex after about half the total runtime of the simulation.

\begin{figure}[htbp]
 \centering
 \caption{Hydromagnetic KH problem  with $Re = 5 \times 10^4, Re_m = 5 \times 10^4, M_A = 15$. Pressure contour lines and velocity streamlines.}
 \label{fig:hmkh}
 \begin{subfigure}{0.24\textwidth} \centering
 \includegraphics[width=\textwidth]{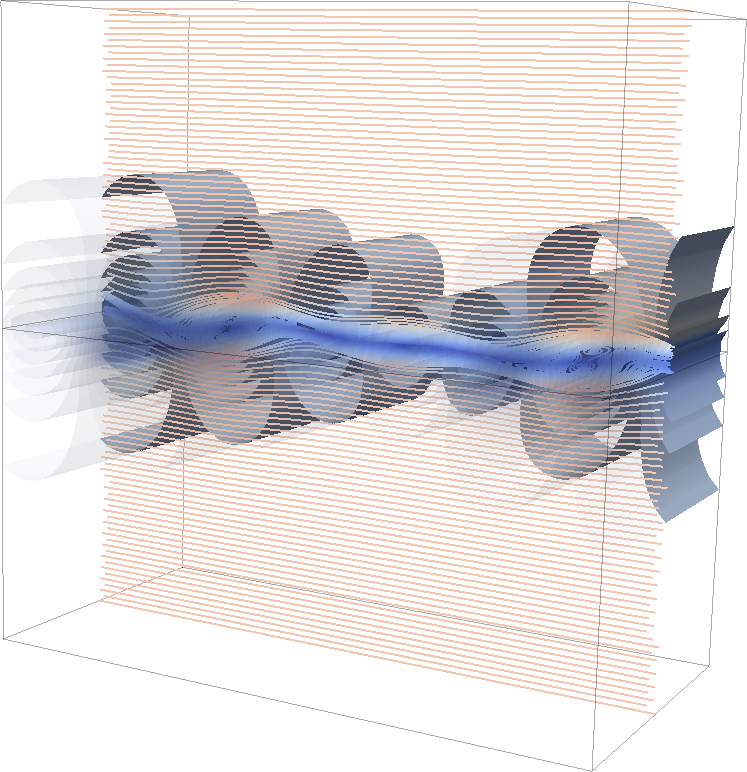}
 \caption{$t=3.0 s$}
 \label{fig:hmkh3}
 \end{subfigure}~%
 \begin{subfigure}{0.24\textwidth} \centering
 \includegraphics[width=\textwidth]{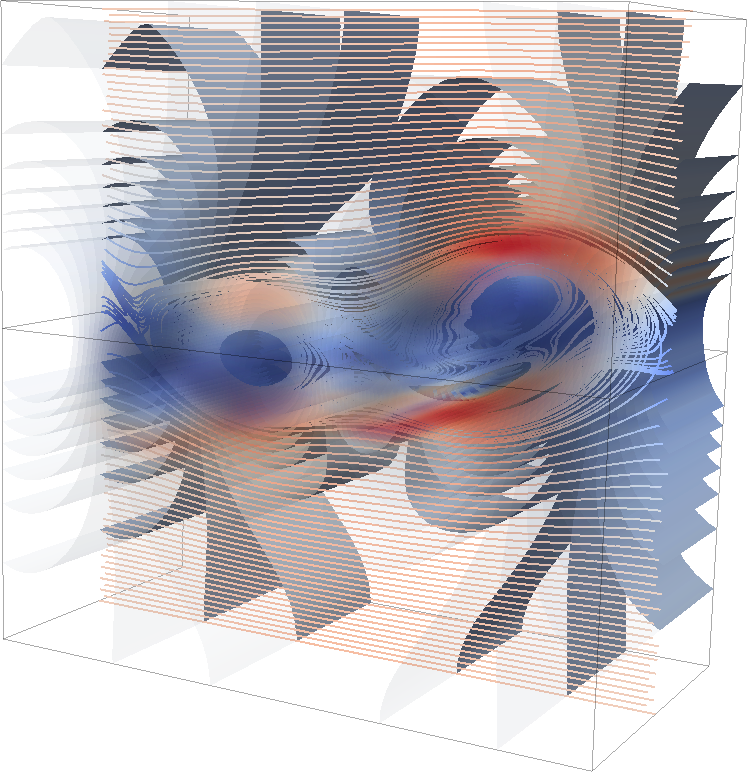}
 \caption{$t=3.5 s$}
 \label{fig:hmkh3x5} 
 \end{subfigure}~%
 \begin{subfigure}{0.24\textwidth} \centering
 \includegraphics[width=\textwidth]{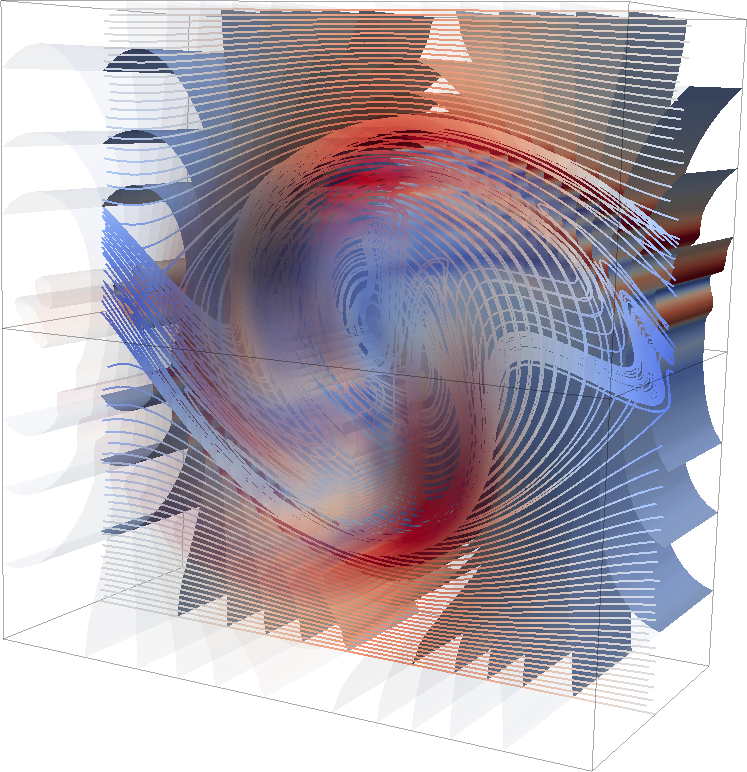}
 \caption{$t=4.0 s$}
 \label{fig:hmkh4} 
 \end{subfigure}~%
 \begin{subfigure}{0.24\textwidth} \centering
 \includegraphics[width=\textwidth]{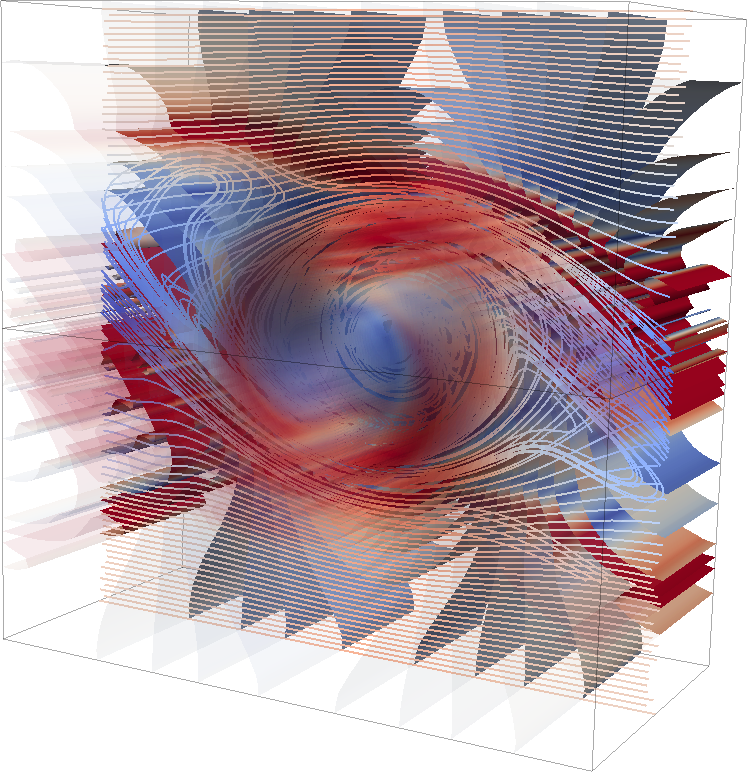}
 \caption{$t=4.5 s$}
 \label{fig:hmkh4x5} 
 \end{subfigure}
\end{figure}


In this numerical study, we perform transient  simulations of the 
problem with  CFL$_{max} = 0.25$ and CFL$_{max} = 0.5$, and study the behavior of the nonlinear- and linear solver. In particular, we compare the number of iterations and timings of the non-restarted GMRES solver using blocked AMG(BGS) as a preconditioner.
The fully-coupled or non-blocked AMG method (FC-AMG) is used as a reference preconditioner.
The linear solver tolerance is set to $10^{-3}$.

Figures \ref{fig:hmkh125_32} and \ref{fig:hmkh25_32} show the solver performance for the different preconditioning strategies over a time sequence for the HMKH problem with a maximum CFL number of $0.25$ and $0.5$ respectively.
We ran the problem both on a $64\times 32\times 16$ mesh on $32$ processors (Intel Broadwell E5-2695 (2.1 GHz), 1 cluster node with $128$ GB RAM) and a $128\times 64\times 32$ mesh on $256$ processors (Intel Broadwell E5-2695 (2.1 GHz), 8 cluster nodes with 128 GB RAM each, Intel Omni-Path high speed interconnect).
Generally, the FC-AMG method needs the least number of iterations,
whereas the blocked AMG(BGS) variant with only 1 coupling iteration needs the highest number of iterations.
Increasing the number of BGS coupling iterations reduces the number of linear iterations getting closer to the reference method.
However, looking at the linear solver time (setup and iteration phase), we see the opposite picture.
The FC-AMG method is the slowest, and the blocked AMG(1 BGS) is the most time efficient method.
Comparing the solver behavior for the different meshes there is a slight increase in the linear iteration count for the finer meshes. So, while the weak scaling is not optimal, the iteration growth with problem size is mild.

\begin{figure}[htbp]
 \centering
 \caption{HMKH example (CFL$_{max} = 0.25$). The left plots show the accumulated linear iterations over time steps. The right plots show the accumulated solution time (setup + iteration phase) per time step.}
 \label{fig:hmkh125_32}
 \begin{subfigure}{1.0\textwidth} \centering
\begin{subfigure}{0.5\textwidth}
\includegraphics[width=0.98\textwidth]{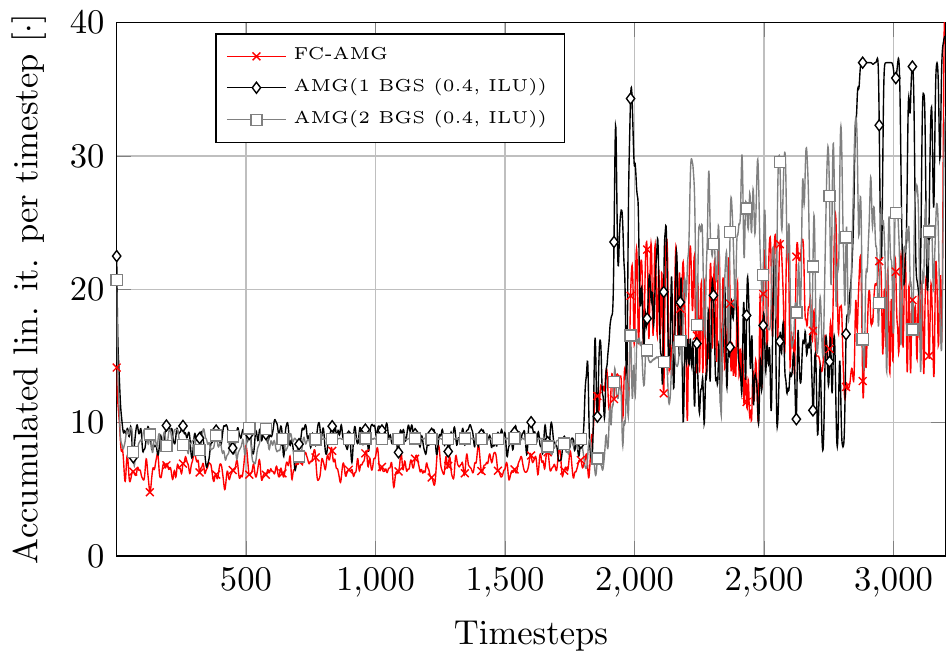}
\end{subfigure}%
\begin{subfigure}{0.5\textwidth}
\includegraphics[width=0.98\textwidth]{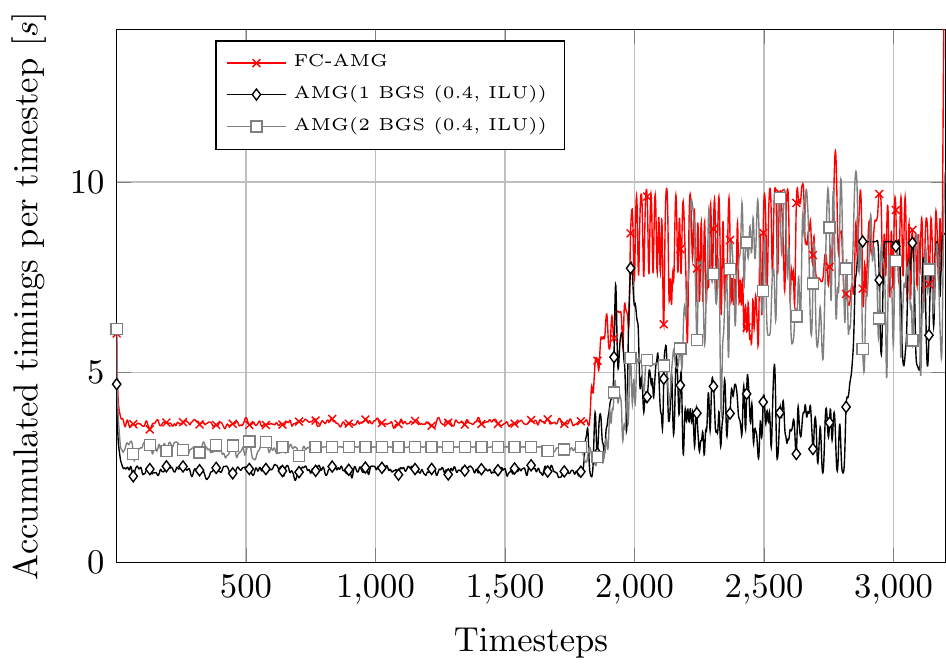}
\end{subfigure}
 \caption{$64\times 32 \times 16$ mesh, $\Delta t = 0.0015625s$, $32$ processors}
\end{subfigure}
\begin{subfigure}{1.0\textwidth} \centering
\begin{subfigure}{0.5\textwidth}
\includegraphics[width=0.98\textwidth]{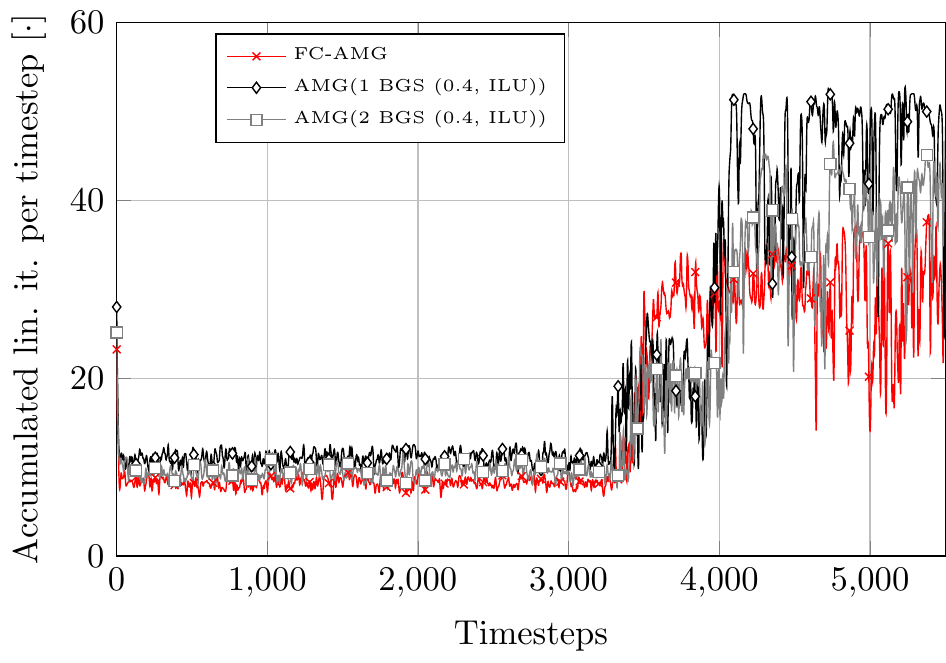}
\end{subfigure}%
\begin{subfigure}{0.5\textwidth}
\includegraphics[width=0.98\textwidth]{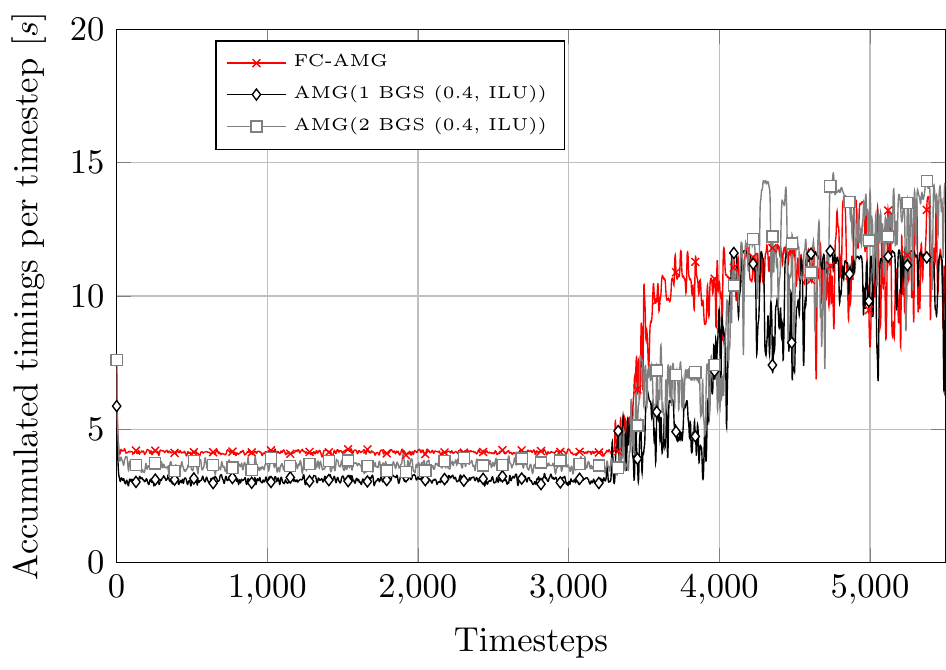}
\end{subfigure}
 \caption{$128\times 64 \times 32$ mesh, $\Delta t = 0.00078125s$, $256$ processors}
\end{subfigure}
\end{figure}

\begin{figure}[htbp]
 \centering
 \caption{HMKH example (CFL$_{max} = 0.5$). The left plots show the accumulated linear iterations over time steps. The right plots show the accumulated solution time (setup + iteration phase) per time step.}
 \label{fig:hmkh25_32}
 \begin{subfigure}{1.0\textwidth} \centering
\begin{subfigure}{0.5\textwidth}
\includegraphics[width=0.98\textwidth]{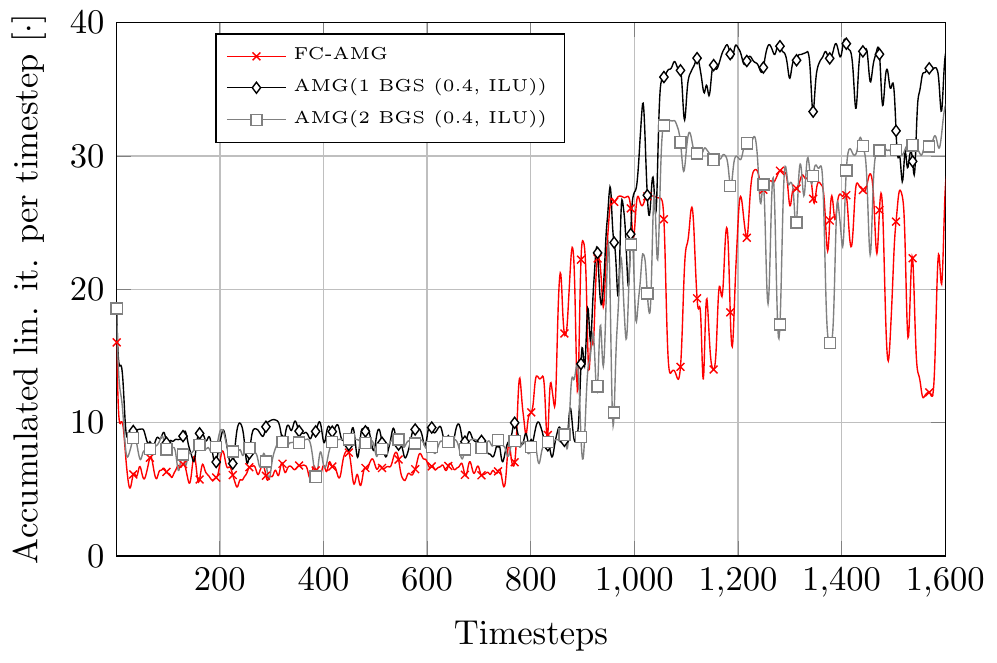}
\end{subfigure}%
\begin{subfigure}{0.5\textwidth}
\includegraphics[width=0.98\textwidth]{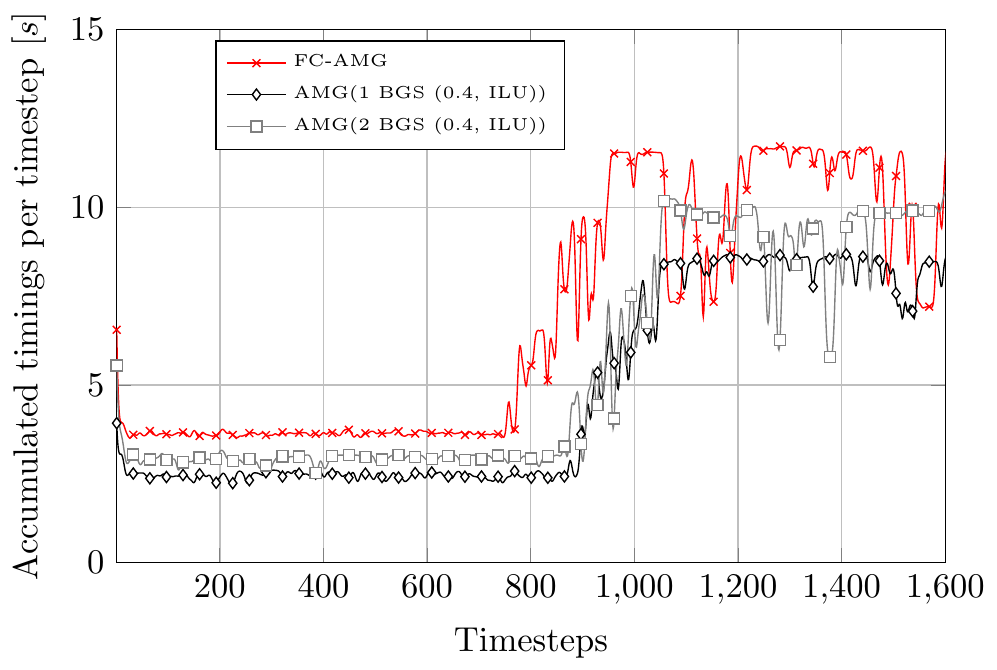}
\end{subfigure}
 \caption{$64\times 32 \times 16$ mesh, $\Delta t = 0.003125s$, $32$ processors}
\end{subfigure}
\begin{subfigure}{1.0\textwidth} \centering
\begin{subfigure}{0.5\textwidth}
\includegraphics[width=0.98\textwidth]{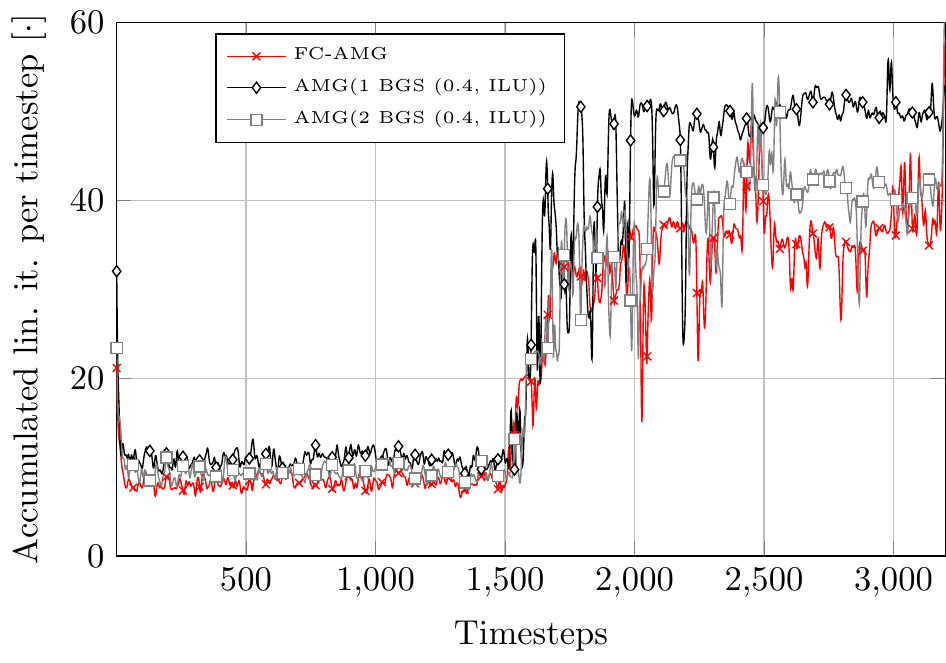}
\end{subfigure}%
\begin{subfigure}{0.5\textwidth}
\includegraphics[width=0.98\textwidth]{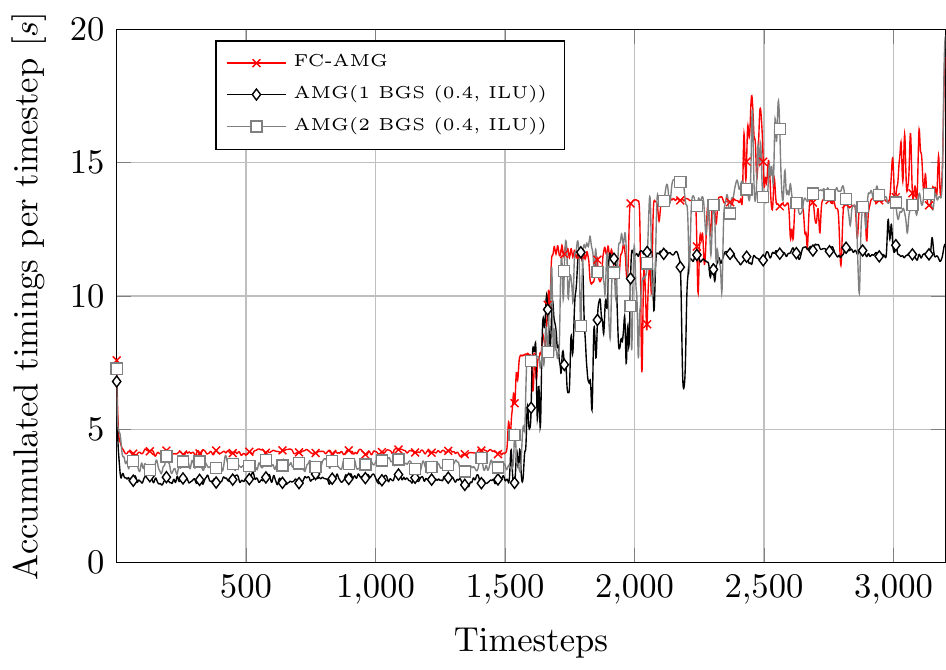}
\end{subfigure}
 \caption{$128\times 64 \times 32$ mesh, $\Delta t = 0.0015625s$, $256$ processors}
\end{subfigure}

\end{figure}

%
%
%
%
%
%

%
%

Numerical results are further summarized in Tables \ref{tab:hmkh125_32} and \ref{tab:hmkh25_32}.
The first column denotes the average number of nonlinear iterations per timestep for the full simulation.
One can see that the number of BGS coupling iterations has some influence on the nonlinear solver, even though there is no clear trend.
The second column represents the average number of linear iterations per nonlinear iteration.
As one would expect, a higher number of BGS coupling iterations reduces the number of linear iterations necessary to solve the problem.
The next three columns give the average setup time, the average solve time and the average overall time for the linear solver per nonlinear iteration.
The AMG(BGS) variants have a clear advantage in the setup costs, but the iteration costs are higher. That is, fully coupled ILU(0) that includes
Navier-Stokes and Maxwell coupling solver provides some convergence benefits, but at the cost of a very large setup time. 
The last three columns show the absolute setup, solve and overall solver time for finishing the simulation.

\begin{table}
\renewcommand{\arraystretch}{1.5}
\caption{HMKH problem for CFL$_{max}=0.25$}
 \label{tab:hmkh125_32}
\resizebox{\columnwidth}{!}{%
\begin{tabular}{p{0.3cm}|l|r|r|rrr|rrr}
& Preconditioner & $\frac{n_N}{n_T}$ & $\frac{n_L}{n_N}$ & $\frac{t_{Se}}{n_N}$ & $\frac{t_{So}}{n_N}$ & $\frac{t_{\Sigma}}{n_N}$ & $t_{Se}$ & $t_{So}$ & $t_{\Sigma}$\\ \hline
\multirow{4}{*}{\small \rotatebox{90}{$64\times 32\times 16$}} 
& FC-AMG & 1.47 & 7.72 & 3.08 & 0.68 & 3.76 & 14443.59 & 3186.92 & 17630.51 \\
& AMG(1 BGS (0.4, ILU)) & 1.37 & 10.18 & 1.46 & 1.11 & 2.57 & 6393.33 & 4854.98 & 11248.30 \\
& AMG(2 BGS (0.4, ILU)) & 1.47 & 9.25 & 1.47 & 1.65 & 3.12 & 6909.57 & 7780.94 & 14690.51 \\
& AMG(3 BGS (0.4, ILU)) & 1.49 & 8.53 & 1.47 & 2.14 & 3.61 & 7034.32 & 10238.93 & 17273.25 \\
\hline
\multirow{4}{*}{\small \rotatebox{90}{$128\times 64\times 32$}} 
& FC-AMG & 1.53 & 10.51 & 3.37 & 0.99 & 4.36 & 28307.94 & 8318.34 & 36626.28 \\
& AMG(1 BGS (0.4, ILU)) & 1.55 & 13.96 & 1.75 & 1.73 & 3.48 & 14884.65 & 14781.09 & 29665.74 \\
& AMG(2 BGS (0.4, ILU)) & 1.51 & 11.87 & 1.74 & 2.39 & 4.13 & 14484.49 & 19942.69 & 34427.18 \\
& AMG(3 BGS (0.4, ILU)) & 1.42 & 10.69 & 1.74 & 3.07 & 4.81 & 13619.80 & 23986.51 & 37606.32 \\
\end{tabular}
}
\tablelegend
\end{table}
\renewcommand{\arraystretch}{1.0}

\begin{table}
\renewcommand{\arraystretch}{1.5}
\caption{HMKH problem for CFL$_{max}=0.5$}
 \label{tab:hmkh25_32}
\resizebox{\columnwidth}{!}{%
\begin{tabular}{p{0.3cm}|l|r|r|rrr|rrr}
& Preconditioner & $\frac{n_N}{n_T}$ & $\frac{n_L}{n_N}$ & $\frac{t_{Se}}{n_N}$ & $\frac{t_{So}}{n_N}$ & $\frac{t_{\Sigma}}{n_N}$ & $t_{Se}$ & $t_{So}$ & $t_{\Sigma}$\\ \hline
\multirow{4}{*}{\small \rotatebox{90}{$64\times 32\times 16$}} 
& FC-AMG & 1.81 & 8.06 & 3.07 & 0.70 & 3.77 & 8878.08 & 2013.02 & 10891.10 \\
& AMG(1 BGS (0.4, ILU)) & 1.77 & 11.19 & 1.47 & 1.24 & 2.71 & 4166.08 & 3513.20 & 7679.28 \\
& AMG(2 BGS (0.4, ILU)) & 1.74 & 9.39 & 1.46 & 1.68 & 3.14 & 4081.49 & 4667.58 & 8749.07 \\
& AMG(3 BGS (0.4, ILU)) & 1.75 & 8.68 & 1.48 & 2.20 & 3.68 & 4143.28 & 6192.65 & 10335.92 \\
\hline
\multirow{4}{*}{\small \rotatebox{90}{$128\times 64\times 32$}} 
& FC-AMG & 1.94 & 11.19 & 3.36 & 1.06 & 4.42 & 20883.40 & 6616.06 & 27499.46 \\
& AMG(1 BGS (0.4, ILU)) & 1.91 & 15.05 & 1.74 & 1.91 & 3.65 & 10691.61 & 11677.32 & 22368.94 \\
& AMG(2 BGS (0.4, ILU)) & 1.89 & 12.74 & 1.74 & 2.58 & 4.32 & 10572.95 & 15626.96 & 26199.91 \\
& AMG(3 BGS (0.4, ILU)) & 1.80 & 11.33 & 1.74 & 3.19 & 4.93 & 10038.69 & 18408.73 & 28447.42 \\
\end{tabular}
}
\tablelegend
\end{table}
\renewcommand{\arraystretch}{1.0}

\subsection{Island coalescence}


\begin{figure}[htbp]
 \centering
 \caption{Structure of the current tubes in 3D island coalescence problem with $\lundquist = 2 \times 10^4$ for the initial condition, and four times in the evolution of the problem $t = 2, 3, 4$. The 3D current tubes have bent in the z-direction and form current sheets.}
 \label{fig:islandpics}
 \begin{subfigure}{0.25\textwidth} \centering
 \includegraphics[width=\textwidth]{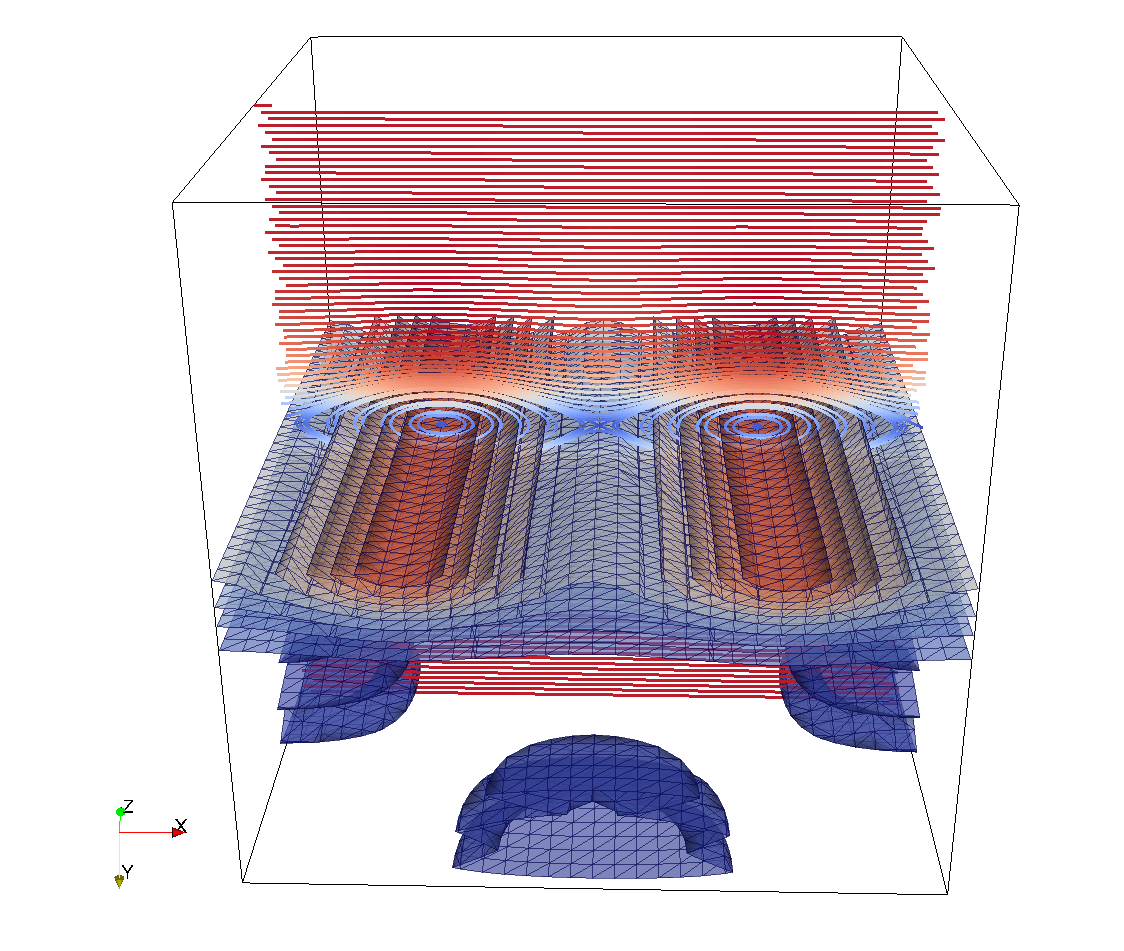}
 \caption{$t=0.0 s$}
 \label{fig:islandpics1}
 \end{subfigure}%
 \begin{subfigure}{0.25\textwidth} \centering
 \includegraphics[width=\textwidth]{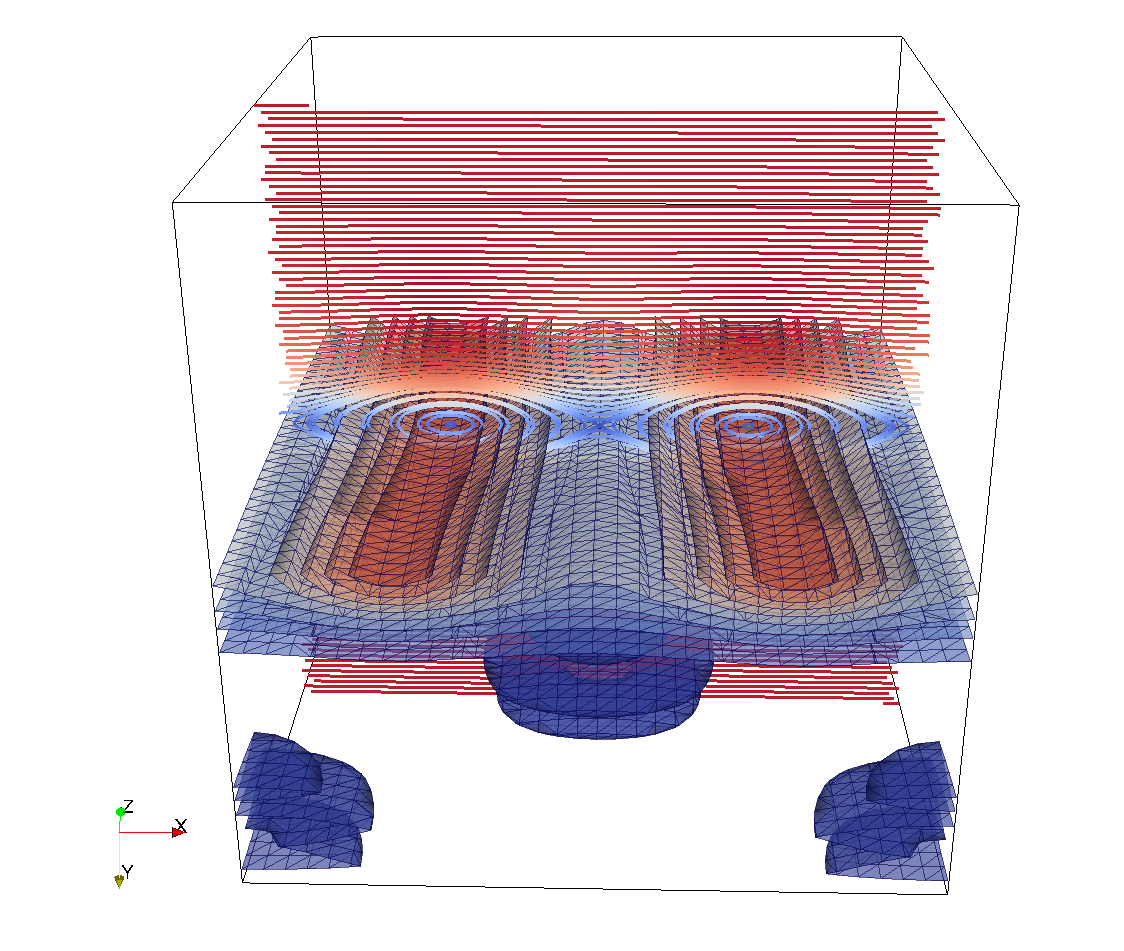}
 \caption{$t=2.0 s$}
 \label{fig:islandpics2} 
 \end{subfigure}%
 \begin{subfigure}{0.25\textwidth} \centering
 \includegraphics[width=\textwidth]{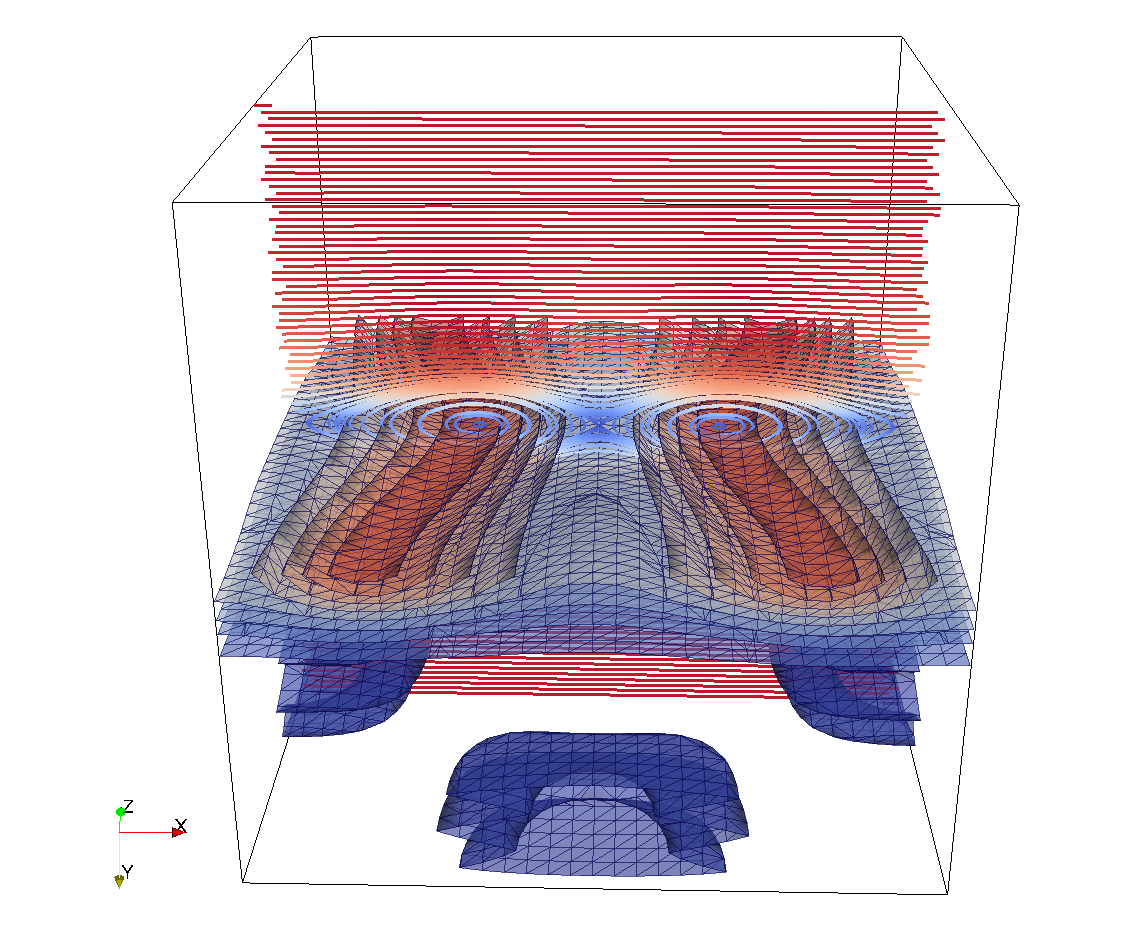}
 \caption{$t=3.0 s$}
 \label{fig:islandpics3} 
 \end{subfigure}%
 \begin{subfigure}{0.25\textwidth} \centering
 \includegraphics[width=\textwidth]{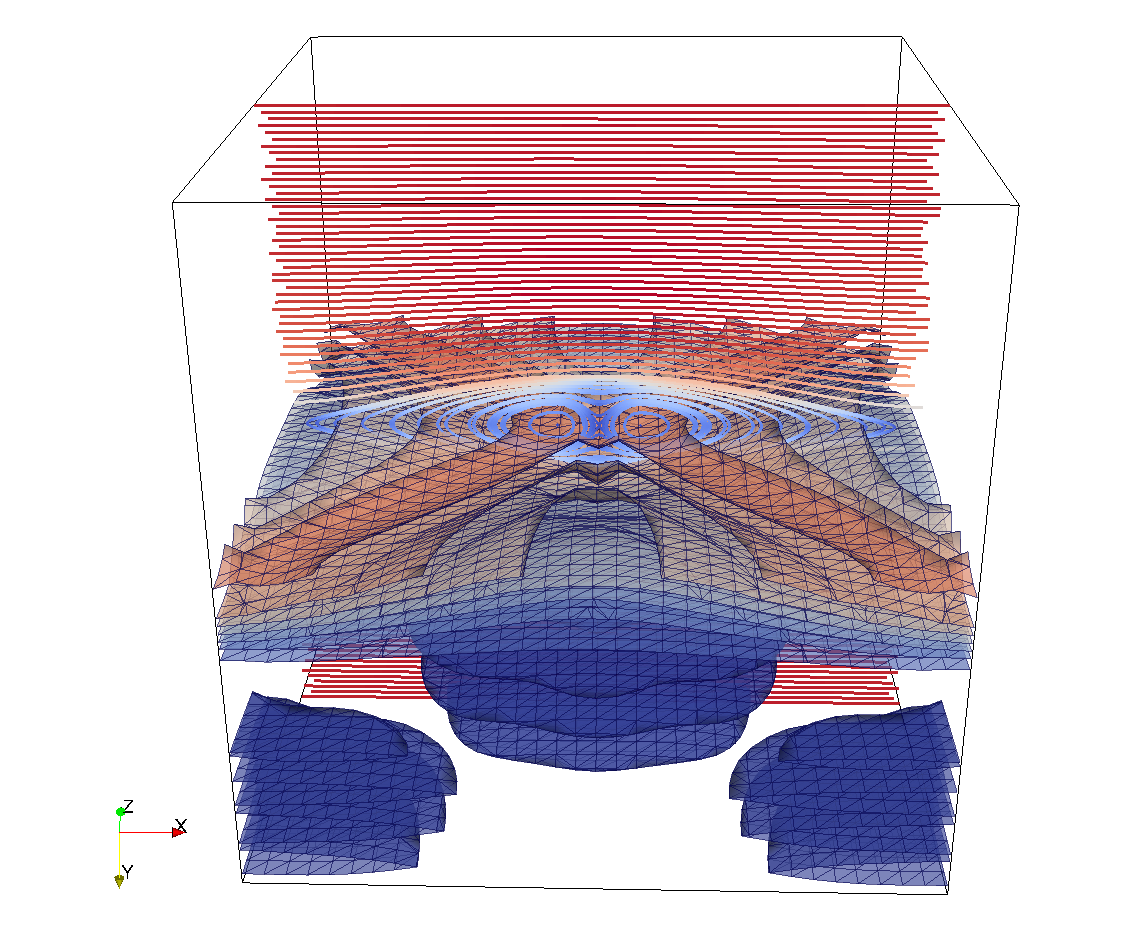}
 \caption{$t=4.0 s$}
 \label{fig:islandpics4} 
 \end{subfigure}
\end{figure}

The island coalescence problem is a prototype problem used to study magnetic reconnection.
Here the initial equilibrium is described by two 3D current tubes (islands in the 2D plane) embedded in
a Harris current sheet (as in the HMKH problem of Section \ref{sec:HMKH})
in a $[-1,1]\times[-1,1]\times[-1,1]$ domain \cite{fadeev-nf-65-ic,chacon-pop-08-3dmhd}.
The initial condition for the island coalescence problem consists of zero fluid velocities ($\vel^0  = {\bf 0}$), 
zero Lagrange multipliers ($\lagr = 0$) and a Fadeev magnetic equilibrium \cite{fadeev-nf-65-ic,chacon-pop-08-3dmhd} that defines the magnetic field $\magn$ and the fluid pressure $\pre$. More details of this setup can be found in \cite{shadid2016scalable}.
The structure of this equilibrium is presented in Figure \ref{fig:islandpics1} with an iso-surface of $\pre$ and iso-lines of $\magn$ at $z=0$. 
The combined magnetic field of the two islands produces Lorentz forces
that pull the islands together.  The dynamics of island coalescence changes as a function of resistivity. For larger resistivities, the x- and o-points 
monotonically approach each other. For low resistivities, fluid-plasma pressure builds up as the islands approach
and a sloshing or bouncing of the o-point position is encountered that leads to lower reconnection rates 
(for more details on the physics see e.g. \cite{biskamp-book-00}).
Figure \ref{fig:islandpics} shows different stages of the reconnection event.
Clearly evident is the formation of the x-point in the intersecting planes between 
the islands (see images at $t = 4$), the development of thin current sheets at that
same x-point location (and the corresponding 3D surface), and the 
movement of the center of the tubes (island o-points) towards the x-point \cite{biskamp-book-00,Knoll-pop-06-ic}.
In this study, we have taken $\dens = 1$, $\visc = \resistivity = 10^{-3}$, $\magnperm = 1$ and,
using the spacing of the o-points, we have $L = 1$, resulting in $\textnormal{Re} = \textnormal{Re}_m = 10^3$.
As in \cite{Knoll-pop-06-ic} these choices imply that the resistivity $\resistivity = 1/\lundquist$,
where $\lundquist$ is the Lundquist number and is defined as $\lundquist = \magnperm L u_A / \resistivity$, 
where $u_A$ is the Alfven velocity. 
We preform transient simulations of the problem with
timestep sizes of $\Delta t \in \{0.05,0.025,0.0125\}$.
The mesh sizes used were $32\times32\times32$, $64\times64\times64$, and $128\times128\times128$
and were run on 8, 64, and 512 processors respectively (Intel Broadwell E5-2695 (2.1 GHz), 2, 16, and 128 cluster nodes with 128 GB RAM each, Intel Omni-Path high speed interconnect).
This provides simulations with a CFL ranging from $1.6$ to $12.8$.
We compare the number of iterations and timings of the non-restarted GMRES solver,
using a relative linear solve tolerance of $\varepsilon=10^{-3}$,
when combined with different preconditioning strategies,
including the fully-coupled or non-blocked AMG method (FC-AMG) as reference and the AMG(BGS) variants.

Figure \ref{fig:island32} shows the solver performance for different preconditioning strategies
over a time sequence for the island coalesce problem with CFL number of $3.2$.
As with the HMKH example, we see that while the AMG(BGS) variants require more iterations than the FC-AMG reference,
the cost per iteration is low enough that the time savings ends up in favor of the AMG(BGS).

Results for the island coalescence problem with various CFL numbers are summarized in Tables~\ref{tab:ic16} through Tables~\ref{tab:ic128}.
The first column denotes the average number of nonlinear iterations per timestep for the full simulation.
The second column represents the average number of linear iterations per nonlinear iteration.
Again, as one would expect, increasing the BGS coupling iterations results in faster convergence or fewer required iterations.
The next three columns give the average setup time,
the average solve time and the average overall time for the linear solver per nonlinear iteration.
The last three columns show the absolute setup, solve and overall solver time for finishing the simulation.
While the FC-AMG boasts faster solve times, it has significant setup costs.
The AMG(BGS) demonstrates a significant reduction in setup time, though
the approach is only slightly faster than the FC-AMG approach
due to the higher solve times.



\begin{figure}[htbp]
 \centering
 \caption{Island coalescence  example (CFL = 3.2). The left plots show the accumulated linear iteratoins over time steps. The right plots show the accumulated solution time (setup + iteration phase) per time step.}
 \label{fig:island32}
 \begin{subfigure}{1.0\textwidth} \centering
\begin{subfigure}{0.5\textwidth}
\includegraphics[width=0.98\textwidth]{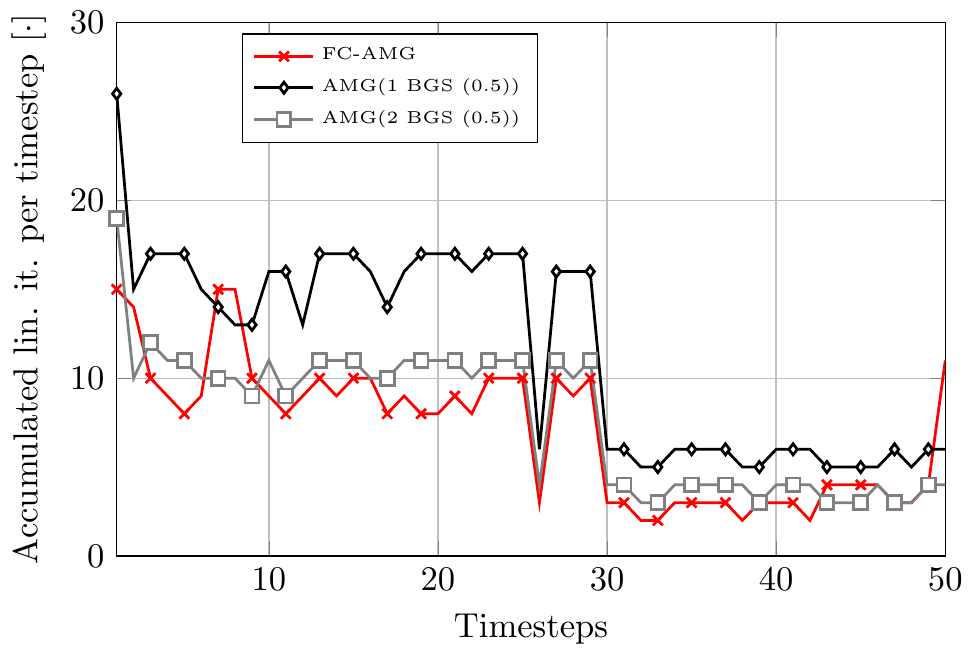}
\end{subfigure}%
\begin{subfigure}{0.5\textwidth}
\includegraphics[width=0.98\textwidth]{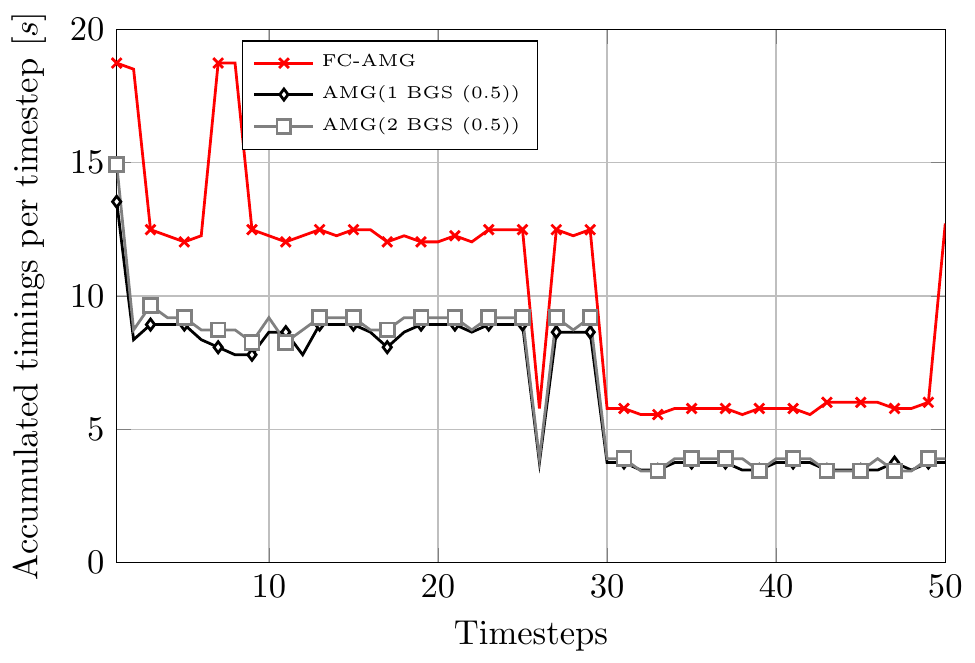}
\end{subfigure}
 \caption{$32\times 32 \times 32$ mesh, $\Delta t = 0.05s$, $8$ processors}
\end{subfigure}
 \begin{subfigure}{1.0\textwidth} \centering
\begin{subfigure}{0.5\textwidth}
\includegraphics[width=0.98\textwidth]{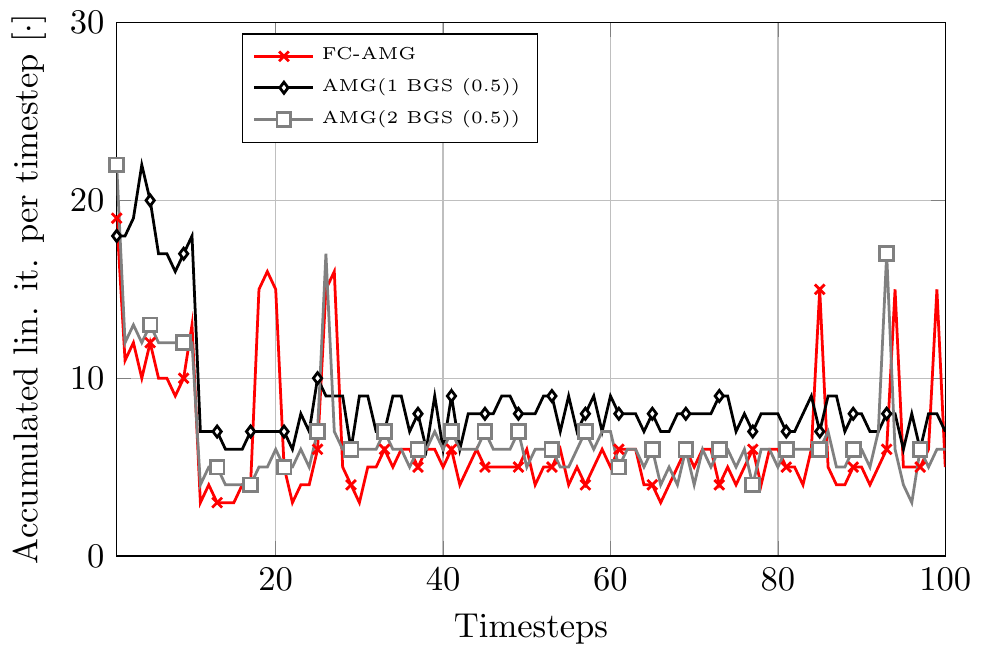}
\end{subfigure}%
\begin{subfigure}{0.5\textwidth}
\includegraphics[width=0.98\textwidth]{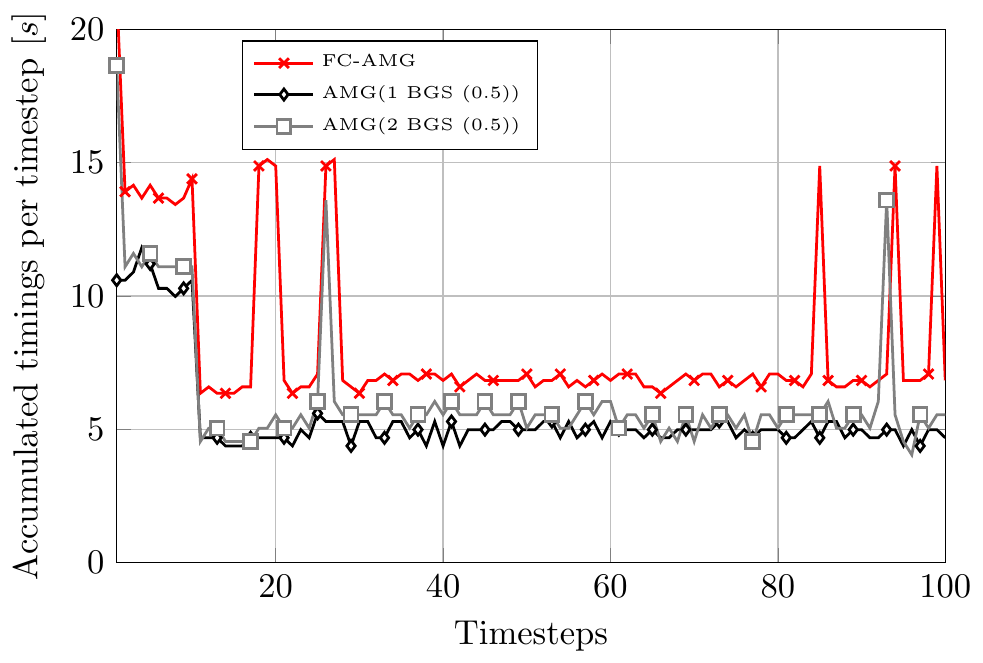}
\end{subfigure}
 \caption{$64\times 64 \times 64$ mesh, $\Delta t = 0.025s$, $64$ processors}
\end{subfigure}
 \begin{subfigure}{1.0\textwidth} \centering
\begin{subfigure}{0.5\textwidth}
\includegraphics[width=0.98\textwidth]{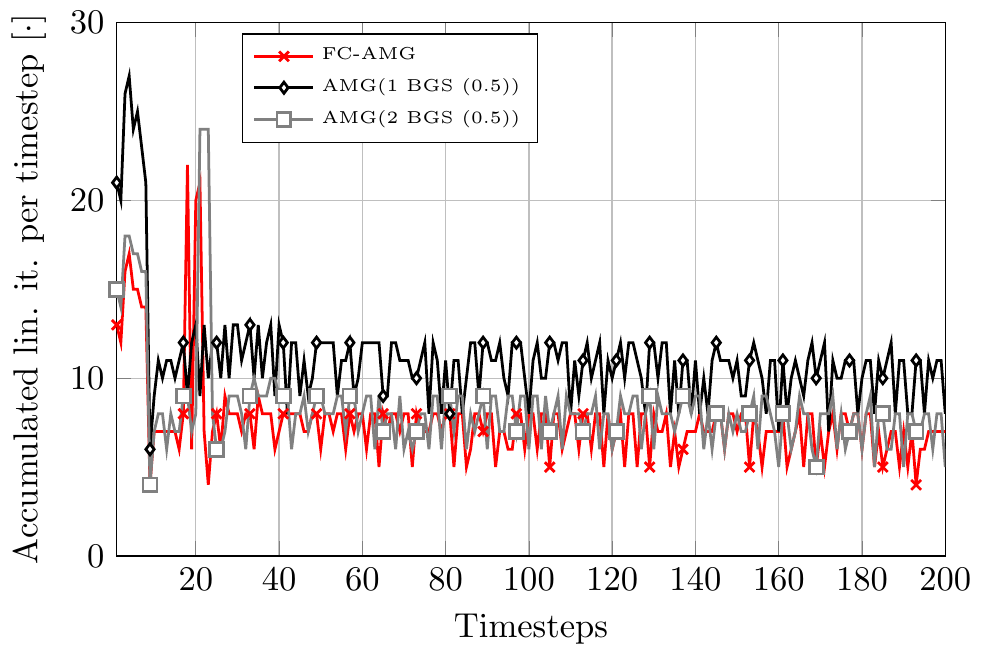}
\end{subfigure}%
\begin{subfigure}{0.5\textwidth}
\includegraphics[width=0.98\textwidth]{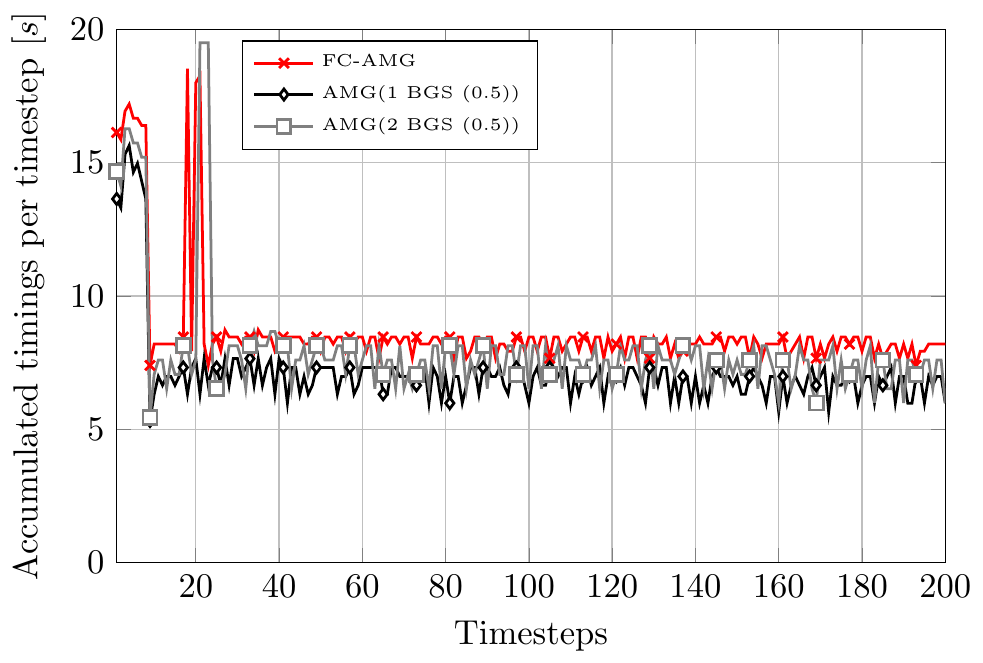}
\end{subfigure}
 \caption{$128\times 128 \times 128$ mesh, $\Delta t = 0.0125s$, $512$ processors}
\end{subfigure}
\end{figure}

\begin{table}
\renewcommand{\arraystretch}{1.5}
\caption{Island coalescence problem for CFL 1.6.}
\label{tab:ic16}
\resizebox{\columnwidth}{!}{%
\begin{tabular}{p{0.3cm}|l|r|r|rrr|rrr}
& Preconditioner & $\frac{n_N}{n_T}$ & $\frac{n_L}{n_N}$ & $\frac{t_{Se}}{n_N}$ & $\frac{t_{So}}{n_N}$ & $\frac{t_{\Sigma}}{n_N}$ & $t_{Se}$ & $t_{So}$ & $t_{\Sigma}$\\ \hline
\multirow{4}{*}{\small \rotatebox{90}{$32\times 32\times 32$}} 
& FC-AMG & 1.16 & 3.96 & 5.07 & 0.91 & 5.98 & 588.12 & 105.67 & 693.79 \\
& AMG(1 BGS (0.4)) & 1.12 & 7.12 & 2.06 & 2.00 & 4.07 & 230.83 & 224.47 & 455.30 \\
& AMG(2 BGS (0.4)) & 1.14 & 5.11 & 2.06 & 2.37 & 4.43 & 234.84 & 269.82 & 504.66 \\
& AMG(3 BGS (0.4)) & 1.15 & 4.29 & 2.06 & 2.74 & 4.80 & 237.36 & 314.98 & 552.34 \\
\hline
\multirow{4}{*}{\small \rotatebox{90}{$64\times 64\times 64$}} 
& FC-AMG & 1.04 & 5.23 & 5.62 & 1.28 & 6.89 & 1174.16 & 266.58 & 1440.74 \\
& AMG(1 BGS (0.4)) & 1.04 & 8.77 & 2.55 & 2.66 & 5.21 & 533.16 & 555.24 & 1088.40 \\
& AMG(2 BGS (0.4)) & 1.04 & 6.64 & 2.55 & 3.31 & 5.86 & 533.16 & 692.22 & 1225.38 \\
& AMG(3 BGS (0.4)) & 1.04 & 5.56 & 2.55 & 3.87 & 6.42 & 532.95 & 808.61 & 1341.56 \\
\end{tabular}
}
\tablelegend
\end{table}
\renewcommand{\arraystretch}{1.0}

\begin{table}
\renewcommand{\arraystretch}{1.5}
\caption{Island coalescence problem for CFL 3.2}
\label{tab:ic32}
\resizebox{\columnwidth}{!}{%
\begin{tabular}{p{0.3cm}|l|r|r|rrr|rrr}
& Preconditioner & $\frac{n_N}{n_T}$ & $\frac{n_L}{n_N}$ & $\frac{t_{Se}}{n_N}$ & $\frac{t_{So}}{n_N}$ & $\frac{t_{\Sigma}}{n_N}$ & $t_{Se}$ & $t_{So}$ & $t_{\Sigma}$\\ \hline
\multirow{4}{*}{\small \rotatebox{90}{$32\times 32\times 32$}} 
& FC-AMG & 1.66 & 4.27 & 5.10 & 0.98 & 6.08 & 423.22 & 81.27 & 504.49 \\
& AMG(1 BGS (0.4)) & 1.60 & 8.54 & 2.07 & 2.42 & 4.48 & 165.44 & 193.33 & 358.77 \\
& AMG(2 BGS (0.4)) & 1.60 & 5.85 & 2.06 & 2.69 & 4.76 & 165.12 & 215.45 & 380.57 \\
& AMG(3 BGS (0.4)) & 1.68 & 4.70 & 2.06 & 2.99 & 5.05 & 173.12 & 251.34 & 424.46 \\
\hline
\multirow{4}{*}{\small \rotatebox{90}{$64\times 64\times 64$}} 
& FC-AMG & 1.19 & 5.35 & 5.63 & 1.29 & 6.92 & 670.33 & 153.54 & 823.87 \\
& AMG(1 BGS (0.4)) & 1.10 & 8.95 & 2.55 & 2.72 & 5.27 & 280.50 & 299.32 & 579.82 \\
& AMG(2 BGS (0.4)) & 1.10 & 6.43 & 2.55 & 3.24 & 5.79 & 280.61 & 356.08 & 636.69 \\
& AMG(3 BGS (0.4)) & 1.13 & 5.36 & 2.56 & 3.71 & 6.27 & 288.72 & 419.54 & 708.26 \\
\hline
\multirow{4}{*}{\small \rotatebox{90}{$128\times 128\times 128$}} 
& FC-AMG & 1.05 & 7.16 & 6.34 & 1.91 & 8.24 & 1336.90 & 402.68 & 1739.57 \\
& AMG(1 BGS (0.4)) & 1.04 & 11.75 & 3.31 & 3.91 & 7.21 & 687.65 & 812.75 & 1500.40 \\
& AMG(2 BGS (0.4)) & 1.05 & 8.73 & 3.32 & 4.69 & 8.01 & 697.83 & 984.75 & 1682.58 \\
& AMG(3 BGS (0.4)) & 1.05 & 7.38 & 3.34 & 5.47 & 8.80 & 700.35 & 1148.35 & 1848.70 \\
\end{tabular}
}
\tablelegend
\end{table}
\renewcommand{\arraystretch}{1.0}

\begin{table}
\renewcommand{\arraystretch}{1.5}
\caption{Island coalescence problem for CFL 6.4}
\label{tab:ic64}
\resizebox{\columnwidth}{!}{%
\begin{tabular}{p{0.3cm}|l|r|r|rrr|rrr}
& Preconditioner & $\frac{n_N}{n_T}$ & $\frac{n_L}{n_N}$ & $\frac{t_{Se}}{n_N}$ & $\frac{t_{So}}{n_N}$ & $\frac{t_{\Sigma}}{n_N}$ & $t_{Se}$ & $t_{So}$ & $t_{\Sigma}$\\ \hline
\multirow{4}{*}{\small \rotatebox{90}{$64\times 64\times 64$}} 
& FC-AMG & 1.50 & 5.61 & 5.60 & 1.35 & 6.96 & 420.15 & 101.61 & 521.76 \\
& AMG(1 BGS (0.4)) & 1.30 & 10.80 & 2.55 & 3.32 & 5.87 & 165.69 & 215.60 & 381.29 \\
& AMG(2 BGS (0.4)) & 1.34 & 6.99 & 2.56 & 3.52 & 6.07 & 171.39 & 235.52 & 406.90 \\
& AMG(3 BGS (0.4)) & 1.54 & 5.86 & 2.56 & 4.07 & 6.63 & 196.81 & 313.53 & 510.34 \\
\hline
\multirow{4}{*}{\small \rotatebox{90}{$128\times 128\times 128$}} 
& FC-AMG & 1.25 & 7.38 & 6.35 & 2.02 & 8.37 & 793.62 & 252.00 & 1045.63 \\
& AMG(1 BGS (0.4)) & 1.09 & 12.27 & 3.31 & 4.09 & 7.40 & 361.01 & 445.69 & 806.70 \\
& AMG(2 BGS (0.4)) & 1.21 & 9.04 & 3.31 & 4.85 & 8.16 & 400.27 & 587.34 & 987.61 \\
& AMG(3 BGS (0.4)) & 1.24 & 7.30 & 3.32 & 5.40 & 8.72 & 411.93 & 669.75 & 1081.67 \\
\end{tabular}
}
\tablelegend
\end{table}
\renewcommand{\arraystretch}{1.0}

\begin{table}
\renewcommand{\arraystretch}{1.5}
\caption{Island coalescence problem for CFL 12.8}
\label{tab:ic128}
\resizebox{\columnwidth}{!}{%
\begin{tabular}{p{0.3cm}|l|r|r|rrr|rrr}
& Preconditioner & $\frac{n_N}{n_T}$ & $\frac{n_L}{n_N}$ & $\frac{t_{Se}}{n_N}$ & $\frac{t_{So}}{n_N}$ & $\frac{t_{\Sigma}}{n_N}$ & $t_{Se}$ & $t_{So}$ & $t_{\Sigma}$\\ \hline
\multirow{4}{*}{\small \rotatebox{90}{$128\times 128\times 128$}} 
& FC-AMG & 1.58 & 7.80 & 6.35 & 2.08 & 8.44 & 502.04 & 164.66 & 666.71 \\
& AMG(1 BGS (0.4)) & 1.34 & 15.48 & 3.31 & 5.16 & 8.47 & 221.84 & 345.86 & 567.69 \\
& AMG(2 BGS (0.4)) & 1.42 & 9.61 & 3.31 & 5.15 & 8.46 & 235.15 & 365.40 & 600.55 \\
& AMG(3 BGS (0.4)) & 1.48 & 7.62 & 3.31 & 5.65 & 8.95 & 244.72 & 417.88 & 662.60 \\
\end{tabular}
}
\tablelegend
\end{table}
\renewcommand{\arraystretch}{1.0}

\subsection{Mixed finite elements}

Next we illustrate a formulation for which the hydrodynamics and electromagnetics systems are
discretized by differing FE spaces. In this example we consider Q2/Q2 VMS
for the hydrodynamics (saddle point  and convective stabilization) and Q1/Q1 VMS for the induction (electromagnetics) systems (saddle point  and convective stabilization). In this problem the difference in order-of-accuracy is motivated by the desire to minimize the overall
computational time while still maintaining higher accuracy for the MHD simulation in appropriate applications.
For example, when a liquid metal is the conducting fluid in an MHD generator,
the flow Reynolds number can be significantly higher than the corresponding magnetic
Reynolds number due to the very high magnetic diffusivity of liquid metals.
In general the low magnetic Reynolds number is indicative of diffusive dominated transport for the magnetics
in the liquid metal. Thus, a mixed discretization with a lower-order approximation for the 
induction subsystem may be appropriate. Other cases that employ disparate discretizations for hydrodynamics and magnetics would be
various forms of structure preserving methods where, for example, nodal FE are employed for flow variables
and face or even edge FE are used for the magnetic field (see e.g. \cite{nedelec:1980,BHRT:2002,schotzau2004mixed,PhillipsShadidCyrElmanPawlowski2016}).

Here, we again consider the MHD Generator problem from Section~\ref{sec:mhdgen} with 
the modest intention of demonstrating the ability of our proposed methods to handle disparate discretizations, which in this case are
Q2/Q2 VMS for the hydrodynamics and Q1/Q1 VMS for the 
induction (electromagnetics) systems.
The system is difficult to approach through standard fully coupled AMG methods due to
the mixed FE spaces with DoFs that are no longer co-located.
The blocked approach outlined in this manuscript allows for the separate construction of
aggregates for the hydrodynamics block and the electromagnetics block.
The monolithic multigrid hierarchy then naturally provides us with coupling between the blocks
on all levels of the hierarchy.

The study is carried out on the same set of physical, geometrical and solver parameters as in Section~\ref{sec:mhdgen},
with varying viscosities $\visc \in \{0.006, 0.01, 0.04 \}$.
The relaxation method is a blocked Gauss-Seidel with a damping parameter of 0.6.
For the sub-block solves, a single iteration of Additive Schwarz with an overlap of 1 was used
to generate approximate sub-block solutions.
The current implementation we are using lacks parallel load rebalancing,
which is problematic 
on higher core counts.
To circumvent this issue, the maximum number of AMG levels was capped at 4 levels,
as further coarsening of the 2048 processor case requires rebalancing.
The coarsest level problem, is still relatively large in the 2048 processor case
(16,384 rows after 3 levels of refinement, or 8 DoFs per processor),
so the coarse level solve is handled with an iteration of the smoother instead of a direct solve.

We explore the weak scaling of the method in Table~\ref{tab:q2q1mhd},
showing iterations and timings for various preconditioner configurations.
For comparison, we also consider an additive Schwarz domain decomposition method (DD-Schwarz)
with overlap 1
for the entire $2\times 2$ block system as a preconditioner.
The only other possible AMG option without using the multiphysics framework
considers the entire system as a scalar PDE. This non-blocked AMG approach
performs so poorly that results are not given. 
The use of Blocked AMG provides a significant reduction in setup time over the
monolithic additive Schwarz, as the block off-diagonal terms are no longer considered in the factorization.
While the number of linear iterations does degrade as the problem size increases,
the ability to apply Blocked AMG to this mixed FE space problem 
provides a significant linear solve time speed up compared to
the use of Additive Schwarz on the entire $2\times2$ block system.
Additional work is needed to better understand smoother and aggregation choices in the mixed FE case.
For example, the increase in iterations as problem size increases for high viscosity indicates potential inefficiencies
with our handling of the Q2 hydrodynamics problem, as the low viscosity case has a more
consistent iteration count as the problem size increases.

\begin{table}
\renewcommand{\arraystretch}{1.4}
\caption{MHD Generator using mixed finite elements.}
\label{tab:q2q1mhd}
\resizebox{\columnwidth}{!}{%
\begin{tabular}{ l| l| l| l| l| l| l| l| l}
&Preconditioner & Processors & visc & $n_N$ & $n_L$  & $n_L/n_N$& $t_{S_e}$ & $t_{So}$ \\ \hline
\multirow{6}{*}{\small \rotatebox{90}{$64\times 32\times 32$}}
&AMG(1 BGS (0.6)) & 32  & 0.04 & 5 & 458& 91.6 & 52.40 & 327.15\\
&AMG(1 BGS (0.6)) & 32  & 0.01 & 5 & 421& 84.2 & 52.27 & 295.11\\
&AMG(1 BGS (0.6)) & 32  & 0.006 & 6 & 970& 161.7 & 62.73 & 709.46\\
&DD-Schwarz               & 32  & 0.04 & 7 & 3016&430.9 & 323.84 & 913.46\\
&DD-Schwarz               & 32  & 0.01 & 11 & 5016&456.0 & 510.68 & 1514.91\\
&DD-Schwarz               & 32  & 0.006 & 21 & 10016& 477.0 & 978.11 & 3029.31\\ \hline
\multirow{6}{*}{\small \rotatebox{90}{$128\times 64\times 64$}}
&AMG(1 BGS (0.6)) & 256  & 0.04 & 5 & 663& 132.6 & 56.97 & 519.58\\
&AMG(1 BGS (0.6)) & 256  & 0.01 & 5 & 475& 95.0 & 56.80 & 357.89\\
&AMG(1 BGS (0.6)) & 256  & 0.006 & 6 & 723& 120.5 & 68.18 & 556.24\\
&DD-Schwarz  & 256  & 0.04 & 15 & 7024& 468.3 & 758.51 & 2200.64\\
&DD-Schwarz  & 256  & 0.01 & 19 & 9024& 474.9 & 963.58 & 2831.34\\
&DD-Schwarz  & 256  & 0.006 & 8 & 3524& 440.5 & 404.60 & 1098.83\\ \hline
\multirow{6}{*}{\small \rotatebox{90}{$256\times 128\times 128$}}
&AMG(1 BGS (0.6)) & 2048  & 0.04 & 5 & 1080& 216.0 & 59.93 & 1003.43\\
&AMG(1 BGS (0.6)) & 2048  & 0.01 & 6 & 993& 165.5 & 71.45 & 856.84\\
&AMG(1 BGS (0.6)) & 2048  & 0.006 & 6 & 944& 157.3& 71.17 & 805.97\\
&DD-Schwarz  & 2048  & 0.04 & 21 & 10032& 477.7& 1103.45 & 3554.11\\
&DD-Schwarz  & 2048  & 0.01 & 22 & 10532& 478.7& 1153.60 & 3729.02\\
&DD-Schwarz  & 2048  & 0.006 & 11 & 5032& 457.4 & 573.59 & 1760.38\\
\end{tabular}
}
\tablelegendmixed
\end{table}

\renewcommand{\arraystretch}{1.0}

Future work will consider correlated coarsening algorithms where the
Q1/Q1 aggregates influence the Q2/Q2 aggregation scheme.
In this example 
there is some partial overlap in the location of DoFs on the mesh.
Some mesh nodes have 8 DoFs, corresponding to four DoFs for hydrodynamics and four DoFs for electromagnetics,
Others only have 4 DoFs, corresponding only to the hydrodynamics.
A natural extension is to force the aggregation scheme 
to  preserve this partial co-location aspect of the hydrodynamics and electromagnetics Dofs.
In future work, we plan to incorporate some ability to partially share aggregation information
between the two AMG invocations. 
In this case, one might share the aggregate root (or central) vertices generated during the 
AMG invocation for electomagnetics. These root vertices could then be used to construct
an initial set of aggregates for the hydrodynamics. As there are more hydrodynamic unknowns, 
many hydrodynamic unknowns might remain unaggregated and so further aggregation 
would be needed to complete the set of aggregates for the hydrodynamics. 

\section{Conclusion}
\label{sec:conclusion}
A new framework for developing multiphysics multigrid preconditioners is
developed and demonstrated on a number of MHD problems.
The key idea is to develop the
multigrid components in a block fashion that mirrors the blocks
in a multiphysics system. Our approach has been to develop block smoothers
and apply them to a multigrid hierarchy constructed using block restriction/prolongation
operators. In many cases, the blocked multiphysics multigrid hierarchy
allows for faster
solution times than a non-blocked approach. For mixed spatial discretizations, the
multiphysics framework provides the only genuine
avenue to leverage pre-existing multigrid software to produce a 
monolithic multigrid preconditioner.  
Here, the AMG engine is invoked multiple times for different sub-blocks and the resulting individual grid 
transfers are combined into one composite operator that
can be employed in a monolithic AMG fashion. 
Ultimately, the run time benefits are much greater because no other multilevel
preconditioning option is available. 
While this paper has focused on specific examples and MHD, the goal of the
framework is to be able to easily construct, adapt, and tailor different
monolithic multigrid preconditioners to various PDE systems.

\bibliographystyle{siamplain}
\bibliography{General_TR_XMHD_ref_prop_FY14,adjoint_stab,bochev,estep,hydro,luis,multilevel,nonlinear,ref,time_int,numerical_MHD,amg_multiphysics} 
\end{document}